\documentclass{amsart}
\usepackage[dvips]{color}
\usepackage[dvips]{graphics}
\usepackage{amsmath}
\usepackage{amsthm}
\usepackage{amsfonts}
\usepackage{amscd}
\usepackage{fancyhdr}
\usepackage{helvet}
\newtheorem{lem}{Lemma}
\newtheorem{theorem}{Theorem}

\newtheorem{cor}{Corollary}
\newtheorem{prp}{Proposition}
\newtheorem{defn}{Definition}

\newtheorem{example}{Example}

\newcommand{\gr}[3]{#1=\left (#2,#3 \right )}

\newcommand{\uned}{\hbox{\kern3pt\raise2.5pt\vbox{\hrule width9pt height 0.3pt}\kern3pt}}
\newcommand{\ded}{\rightarrow}
\newcommand{\bded}{\leftrightarrow}

\newcommand{\ind}{\perp\!\!\!\perp}
\newcommand{\cind}[3]{#1 \ind #2|#3}

\newcommand{\path}[2]{\raisebox{-.6ex}{{\footnotesize $#1$}}{\Large
$\pi$}\raisebox{-.6ex}{{\footnotesize $#2$ }}}
\newcommand{\epath}[2]{\raisebox{-.6ex}{\footnotesize $#1$}
{\text{\Large $\pi$}}\raisebox{-.6ex}{\footnotesize $#2$}}

\newcommand{\cm}[3]{#1\Bigl[^{#2}_{#3}}

\newcommand{\rhsq}[3]{\rho^2_{#1#2\mid #3}}

\newcommand{\indrels}[1]{\mathfrak{I}\left(#1\right)}
\newcommand{\trip}[3]{\ensuremath{\langle #1,
#2\mid #3\rangle}} 

\usepackage{array}

\usepackage[sort]{natbib}
\usepackage{fullpage}
\usepackage{pstricks,pst-node}
\usepackage{subfigure}
\usepackage{graphicx}
\newcommand{\cn}[1]{\circlenode*{.05cm}{#1}}

\begin{document}

\title[]{Qualitative inequalities for squared partial correlations of a Gaussian random vector}

\author{Sanjay Chaudhuri}



\address{Department of Statistics and Applied probability, National University of Singapore, Singapore, $117546$.}
\email{sanjay@stat.nus.edu.sg}

\thanks{This research was partially supported by Grant R-155-000-081-112 from National University of Singapore.}



\date{\today}




\begin{abstract}
We describe various sets of conditional independence relationships, sufficient for qualitatively comparing non-vanishing squared partial correlations of a Gaussian random vector. These sufficient conditions are satisfied by several graphical Markov models.  
Rules for comparing degree of association among the vertices of such Gaussian graphical models are also developed.
We apply these rules to compare conditional dependencies on Gaussian trees. 
In particular for trees, we show that such dependence can be completely characterised by the length of the paths joining the dependent vertices to each other and to the vertices conditioned on.  We also apply our results to postulate rules for model selection for polytree models.  
Our rules apply to mutual information of Gaussian random vectors as well.  
      





\end{abstract}

\keywords{Inequalities, graphical Markov models, mutual information, squared partial correlation, tree models.}

\maketitle

\section{Introduction}

In graphical Markov models literature, several attempts have been made to characterise the degree of conditional association among the vertices by the structure of the underlying graph. 
Such knowledge is considered useful in model selection. 
For example, \citet{cheng1} describe an algorithm of model selection for directed acyclic graphs (DAG) which assumes that the mutual information has a monotone relationship with certain structure based length of the path.  
Examples \citep{meek1} show that such a \emph{monotone DAG faithfulness} property or a similar \emph{compound monotone DAG faithfulness} property do not hold even for simple binary DAGs. 
In fact, except in some specific cases e.g. \citet{greenland1} in epidemiology, \citet[causal pipes]{spirtes:etal:2000} in causal analysis, no result is known in this context.   

A more general problem is to order the squared partial correlation coefficients among the components of a Gaussian random vector.  For these random vectors, squared partial correlation coefficients completely measure the degree of association between its components conditional on a subset of the components. 
This measure is a polynomial in the entries of their covariance matrices. Thus in many situations it is beneficial to be able to order squared partial correlation coefficients in a way, such that the ordering does not depend on the specific values of the covariances.  

Simple counter-examples show that such \emph{qualitative} comparisons cannot hold unless the covariance matrix belongs to certain subsets of positive definite matrices.  In this article, we specify such subsets by conditional independence relationships. 
For a graphical Markov model validity of such relationships can be simply read off from the underlying graph. Thus rules for comparing degree of association on various Gaussian graphical models can be developed.   
 
In this article we show that, certain conditional independence relationships holding, suitable squared partial correlations can be qualitatively compared.  We make two kinds of comparisons. In the first, the set of components conditioned on (conditionate) are kept fixed and we change the dependent vertices (correlates).  
More importantly, in the second, we fix the two correlates and compare their degree of dependence by varying the conditionates. The sufficient conditional independence relationships are satisfied by several graphical Markov models.
 Using relevant \emph{separation} criteria (e.g. separation for undirected graphs (UG) (see Definition \ref{defn:sep}), d-separation for DAGs \citep{verma_pearl_1990} (see Definition \ref{defn:dconn}), m-separation for mixed ancestral graphs (MAGs) (see supplement) \citep{thomas1} etc., 
we postulate sufficient structural conditions for comparing conditional association on them.
We emphasize that the specific graphical Markov models are used as illustrations.  Our results apply to a much wider class of models. 
Furthermore, using the fact that for tree and polytree (DAGs without any undirected cycles either or singly connected directed acyclic graphs) models, any two connected components have exactly one path joining them, these structural criteria can be simplified to path based rules for comparison. 
We discuss such rules for trees in details, where it is also shown that our rules for comparing the squared partial correlations are complete. 

The inequalities discussed here have theoretical interest as new properties of Gaussian random vectors and directly translate to corresponding conditional non-Shannon type information inequalities \citep{zhang:yeung:1997,matus:2006,matus:2007}.  
\citet{matus:2005} considers implications of one set of conditional independence relations on other conditional independencies for Gaussian random vectors. Furthermore, he describes a way to determine such implications using the ring of polynomials generated by the entries of the correlation matrices with some additional indeterminates. Our results describe some polynomial inequalities these rings satisfy.   

Our main motivation comes from the Gaussian graphical Markov models. 
These results are canonical and sufficient to postulate structure based rules to order dependencies on several of them.   
We improve upon \citet{sctsr3,scphd}, who only consider polytree models. 
These results can be used in determining the distortion effects \citep{wermuth_cox_2008} and monotonic effects \citep{vanderweele_robins_2007,vanderweele_robins} of confounded variables in epidemiology and causal network analysis (see also \citet{greenland_pearl_2011}). 
We postulate necessary and sufficient conditions for determining structures on a class of polytree models.  These conditions can be directly applied in model selection, specially in mapping river flow and drainage networks where such polytree models occur naturally \citep{rodriguez-iturbe_rinaldo_book}.  
In real data analysis, these inequalities would be useful for model selection, specially among various graphical Markov models \citep{cheng1,shimizu:hoyer:etal:2006}.  For these models our results would translate to hypothesis connected to the structure of the graph. These hypothesis can be tested from the observed data. 
 Structure based inequalities may also be used as constraints in estimation with missing values.  They are also relevant in choosing prior distributions in Bayesian procedures. 
The qualitative bounds can be used in selecting stratifying variables in designing surveys, gathering most relevant information in forensic sciences and building strategies for constrained searches. 
Further, these results may have applications in designing effective updating and blocking strategies in Gibbs sampling and Markov chain monte carlo procedures (see eg. \citet{roberts_sahu_1997} etc).

~\hfill\\
\section{Squared partial correlation inequalities}\label{sec:canon}
Suppose $V\sim N\left(\mu,\Sigma\right)$ with a positive definite $\Sigma$.  Let $a$, $b$, $c$, $c^{\prime}$, $z$ , $z^{\prime}$, $x$ etc. be the components and $B$, $Z$ etc. be the subsets of components of $V$.  In this article $V$ will also denote the vertex set of the underlying graph (see supplement for more details). Let $\emptyset$ denote the empty set.

The squared partial correlation coefficient ($\rho^2_{ac|Z}$) between $a$ and $c$ conditional on $Z$ is defined by:

\begin{equation}\label{eq:rho}
\rhsq{a}{c}{Z}=\frac{\left(\sigma_{ac}-\Sigma_{aZ}\Sigma_{ZZ}^{-1}\sigma_{cZ}\right)^2}{\left(\sigma_{aa}-\Sigma_{aZ}\Sigma_{ZZ}^{-1}\sigma_{aZ}\right)\left(\sigma_{cc}-\Sigma_{cZ}\Sigma_{ZZ}^{-1}\sigma_{cZ}\right)}=1-e^{-2Inf\left(\cind{a}{c}{Z}\right)}.
\end{equation}

Here $\sigma_{ab}$ and $\Sigma_{aZ}$ respectively denote the $(a,b)$th element and $a\times Z$ submatrix of $\Sigma$.  $Inf\left(\cind{a}{c}{Z}\right)$ is the mutual information \citep[information proper]{whit} of $a$ and $c$ given $Z$.  From \eqref{eq:rho} it follows that  the mutual information is a monotone increasing function of the corresponding squared partial correlation.  
Thus the qualitative inequalities for $\rhsq{a}{c}{Z}$ presented below applies to $Inf\left(\cind{a}{c}{Z}\right)$ as well.

\subsection{Comparing conditional dependence with a fixed conditionate}\label{sec:comp1}
We first fix a subset $Z$ to be conditioned and one correlate $a$. The squared partial correlation is compared by changing the other correlate from $c$ to $c^{\prime}$.

\begin{theorem}\label{lem:maincomp1}
Suppose $\cind{c^{\prime}}{a}{cZ}$, then $\rhsq{a}{c^{\prime}}{Z}\le\rhsq{a}{c}{Z}$.
\end{theorem}

Theorem \ref{lem:maincomp1} is a conditional version of the well-known \emph{information inequality} \citep{cover} and holds in general for mutual information of any distribution.
  For graphical Markov models the condition holds if $c^{\prime}$ is separated from $a$ given $c$ and $Z$. Further, for trees the condition is satisfied if $c$ lies on the path joining $a$ and $c^{\prime}$. 
Thus longer path implies weaker dependence in this case.  

\begin{figure}[t]
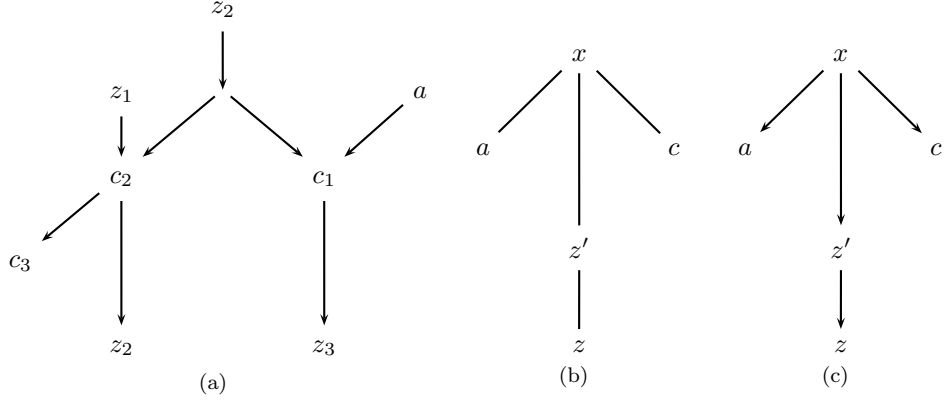

\begin{center}
\subfigure[\label{fig:fxdcond}]{
\input{fxdcond.tex}
}\hspace{1pt}
\subfigure[\label{fig:figcn1}]{
\input{figcn1.tex}
}
\hspace{1pt}
\subfigure[\label{fig:figap1b2}]{
\input{figap1b2.tex}
}
\caption{\ref{fig:fxdcond} A polytree, $Z=\{z_1,z_2,z_3,z_4\}$, $\rhsq{a}{c_3}{Z}\le\rhsq{a}{c_2}{Z}$, however $\rhsq{a}{c_3}{Z}\le\rhsq{a}{c_1}{Z}$ may not hold.  
\ref{fig:figcn1} an UG satisfying the conditions of Theorem \ref{ap:1b}, $\rho^2_{ac}\ge\rhsq{a}{c}{z}\ge\rhsq{a}{c}{z^{\prime}}$ and $\rho^2_{ax}\ge\rhsq{a}{x}{z}\ge\rhsq{a}{x}{z^{\prime}}$. 
Further, from Theorem \ref{lem:maincomp1}, $\rho^2_{ac}\le\rho^2_{ax}$, $\rhsq{a}{c}{z}\le\rhsq{a}{x}{z}$ and $\rhsq{a}{c}{z^{\prime}}\le\rhsq{a}{x}{z^{\prime}}$. Exactly the same conclusions hold on the DAG in \ref{fig:figap1b2}.}
\end{center}
\end{figure}

For polytree models the condition depends on the arrangement of the arrows on the path joining $a$, $c$ and $c^{\prime}$. The condition is satisfied if two arrowheads do not meet at $c$ on the path joining $a$ and $c^{\prime}$, 
(ie. $c$ is not a \emph{collider} on the path joining $a$ and $c^{\prime}$, see Definition \ref{defn:coll}).  As for example, in Figure \ref{fig:fxdcond} with $Z=\left\{z_1,z_2,z_3,z_4\right\}$, using the d-separation criterion (see Definition \ref{defn:dconn}) we get, $\cind{c_3}{a}{Zc_2}$.  
Theorem \ref{lem:maincomp1} ensures that $\rhsq{a}{c_3}{Z}\le\rhsq{a}{c_2}{Z}$. The same d-separation criterion however implies that $c_3\not\ind a|Zc_1$, 
so there is no guaranty the $\rhsq{a}{c_1}{Z}$ would be larger than $\rhsq{a}{c_3}{Z}$.  This partially justifies the intuitive argument given in \citet{greenland1} (see also \citet{greenland_pearl_2011}).

\subsection{Comparing conditional dependence with fixed correlates}\label{sec:appendix} Here two components $a$ and $c$ of $V$ are held fixed.  We consider the variation in $\rhsq{a}{c}{Z}$ for different subsets $Z$ of $V$.  
  Depending on the nature of pairwise unconditional association between $a$, $c$ and the sets conditioned on, three situations may arise.


\begin{figure}[t]
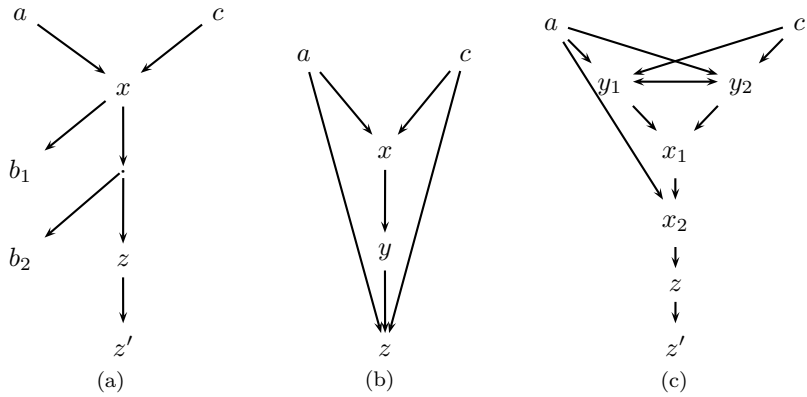

\begin{center}
\subfigure[\label{fig:ap2b}]{
\input{figap2b.tex}
}\hspace{10pt}
\subfigure[\label{nanny}]{
\input{nanny.tex}
}\hspace{10pt}
\subfigure[\label{mag}]{
\input{mag.tex}
}
\end{center}
\caption{Graphical models satisfying the conditions of Theorem \ref{ap:2b}. In each graph $a\ind c$. The graph in \ref{fig:ap2b} is a polytree. Here $B=\{b_1,b_2\}$ and $\rhsq{a}{c}{B}\le\rhsq{a}{c}{Bz^{\prime}}\le\rhsq{a}{c}{Bz}$ holds.  In \ref{nanny} it follows that $\rhsq{a}{c}{y}\le\rhsq{a}{c}{x}$ (cf. \citet{wermuth_cox_2008}). The graph in \ref{mag} is a mixed ancestral graph \citep{thomas1} where $\rhsq{a}{c}{z^{\prime}}\le\rhsq{a}{c}{z}$ always holds.}
\label{fig:ap3ba}
\end{figure}
    
\subsubsection{\bf Situation $1$.} The components $a$, $c$, $z$ and $z^{\prime}$ are unconditionally pairwise dependent. 
\begin{theorem}\label{ap:1b}
Suppose for some $x$, $\cind{a}{c}{x}$ and $\cind{ac}{z}{x}$.  Then $\rhsq{a}{c}{z}\le\rho^2_{ac}$. In addition, if $\cind{ac}{z^{\prime}}{z}$, then $\rhsq{a}{c}{z}\le\rhsq{a}{c}{z^{\prime}}\le\rho^2_{ac}$.

\end{theorem}

The conditions of Theorem \ref{ap:1b} can be represented by several graphical Markov models, eg. undirected graphs, directed acyclic graphs etc.  The conditional independence conditions imply that $a$, $c$ and $z$ have to be pairwise separated given $x$ and $z^{\prime}$ has to be separated from $a$ and $c$ given $z$. 


The first part shows that under these conditions the dependence of $a$ on $c$ always reduces on conditioning.  For tree and polytree models the conclusion of the second part can be intuitively explained. 
Notice that, by assumption $\rho^2_{ac}\ge\rhsq{a}{c}{x}=0$ and the separation criteria imply that $z^{\prime}$ is farther away from $x$ than $z$. Thus $z^{\prime}$ has less information about $x$ than $z$. So $\rhsq{a}{c}{z^{\prime}}$ should be closer to $\rho^2_{ac}$ than $\rhsq{a}{c}{z}$.
In other words, conditioning on the vertices farther away from the path between $a$ and $c$ increases the degree of association.

\subsubsection{\bf Situation $2$.} The correlates $a$ and $c$ are independent, but both are dependent on the sets conditioned on.  
\begin{theorem}\label{ap:2b}
Suppose $a\ind c$ and for some $x$, the condition $\cind{ac}{zB}{x}$ holds.  Then $\rhsq{a}{c}{B}\le\rhsq{a}{c}{Bz}$.  Moreover, if $\cind{z^{\prime}}{acB}{z}$ holds, then $\rhsq{a}{c}{B}\le\rhsq{a}{c}{Bz^{\prime}}\le\rhsq{a}{c}{Bz}$.
\end{theorem}

By assumption $0=\rho^2_{ac}\le\rhsq{a}{c}{B}$.  Thus the first conclusion implies that conditioning on a larger set implies stronger association.
 On an UG, the condition $a\ind c$ implies that $a$ and $c$ cannot be connected.  Thus UGs are not useful to represent the conditions in Theorem \ref{ap:2b}.
  They are satisfied by several other graphical Markov models like DAGs, MAGs etc.

For polytree models (See Figure \ref{fig:ap2b}) the conclusions of Theorem \ref{ap:2b} can be intuitively explained as well. As before, one can conclude $z^{\prime}$ is farther away from $x$ and therefore has less information about $x$ than $z$, 
$\rhsq{a}{c}{x}\ne 0$ but $\rho^2_{ac}=0$. Thus by the same argument as for Theorem \ref{ap:1b}, conditioning on $B$ and $z^{\prime}$ should produce weaker association than $B$ and $z$.    

In the graph in Figure \ref{nanny} the marginal covariance matrix of $a$, $c$, $x$ and $y$ satisfy the conditions of Theorem \ref{ap:2b}. Thus, $\rhsq{a}{c}{y}\le\rhsq{a}{c}{x}$.  The graph in Figure \ref{mag} is a mixed ancestral graph (notice the $\leftrightarrow$ edge between $y_1$ and $y_2$ \citep{thomas1}). 
Here the marginal covariance matrix of $a$, $c$, $x_2$, $z$ and $z^{\prime}$ would satisfy the conditions of Theorem \ref{ap:2b} (see Appendix \ref{sec:mags}). So we conclude that $\rhsq{a}{c}{z^{\prime}}\le\rhsq{a}{c}{z}\le\rhsq{a}{c}{x_2}$.

 




\begin{figure}[t]
\begin{center}
\subfigure[\label{fig:ap3}]{
\input{figap3.tex}
}\hspace{1pt}
\subfigure[\label{fig:ap3b1ca}]{
\input{figap3b1ca.tex}
}\hspace{1pt}
\subfigure[\label{fig:ap3b1c}]{
\input{figap3b1c.tex}
}
\end{center}
\caption{Graphical models satisfying the conditions of Theorem \ref{ap:3b}. Each model satisfies the condition $(i)$ of the theorem. \ref{fig:ap3} is a polytree on which $\cind{aczz^{\prime}}{\{b_1,b_2\}}{x}$ holds. In \ref{fig:ap3b1ca}, $\cind{ac}{b}{x}$, but $ac\not\ind b|zx$. In \ref{fig:ap3b1c}, $\cind{ac}{b}{zx}$ but $ac\not\ind b|x$. From Theorem \ref{ap:3b} it follows that $\rhsq{a}{c}{B}\le\rhsq{a}{c}{Bz^{\prime}}\le\rhsq{a}{c}{Bz}$.}
\label{fig:ap3bb}
\end{figure} 

\subsubsection{\bf Situation $3$.} At least one of $a$ and $c$ is independent of both the sets conditioned on. 


\begin{theorem}\label{ap:3b}
Suppose $a\ind z$. Let for some $x$, $\Sigma$ satisfies one of the following two $( (i), (ii) )$ conditions: 
\begin{list}{}{}
\item[$(i)$] $c\ind az$ and one of the following six conditions $(a)$ $\cind{az}{B}{x}$, $(b)$ $\cind{az}{B}{cx}$, $(c)$ $\cind{cz}{B}{x}$, $(d)$ $\cind{cz}{B}{ax}$, $(e)$ $\cind{ac}{B}{x}$ and $(f)$ $\cind{ac}{B}{xz}$ holds, 
\item[$(ii)$] $\cind{az}{cB}{x}$.
\end{list}
Then $\rhsq{a}{c}{B}\le\rhsq{a}{c}{Bz}$.  Further, if $\cind{z^{\prime}}{acB}{z}$ holds, then in both cases, $\rhsq{a}{c}{B}\le\rhsq{a}{c}{Bz^{\prime}}\le\rhsq{a}{c}{Bz}$. 
\end{theorem} 
\begin{figure}[t]
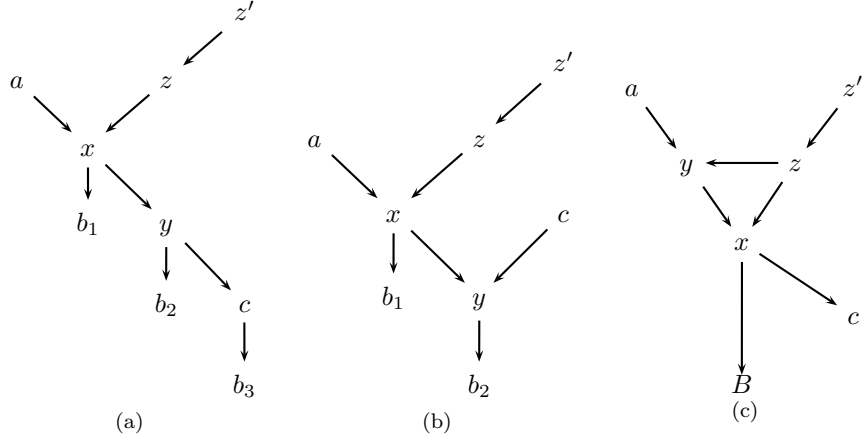

\subfigure[\label{fig:ap3b2}]{
\input{figap3b2.tex}
}\hspace{1pt}
\subfigure[\label{fig:ap3b1}]{
\input{figap3b1.tex}
}\hspace{1pt}
\subfigure[\label{marg}]{
\input{marg.tex}
}
\caption{Graphical models satisfying the conditions $(ii)$ of Theorem \ref{ap:3b}. In each graph the condition $\cind{az}{cB}{x}$ holds. In Figure \ref{fig:ap3b1} $\cind{az}{\{b_1,b_2\}}{cx}$ also holds.  The graphs in \ref{fig:ap3b2} ($B=\{b_1,b_2,b_3\}$) and \ref{fig:ap3b1} ($B=\{b_2,b_3\}$) are polytrees. On each $\rhsq{a}{c}{B}\le\rhsq{a}{c}{Bz^{\prime}}\le\rhsq{a}{c}{Bz}$ hold.}
\label{fig:ap3b}
\end{figure} 
The difference between the conditions $(i)$ and $(ii)$ in Theorem \ref{ap:3b} is illustrated in Figure \ref{fig:ap3} and \ref{fig:ap3b2}.  Under condition $(i)$, $c\ind z$ but the relation $c\not\ind z|x$ does not necessarily hold .  On the other hand, under condition $(ii)$, $\cind{c}{z}{x}$ but $c$ may not be independent $z$ unconditionally.   

The six conditions in $(i)$ are in general distinct. As for example, from m-connection rules \citep{thomas1} the MAG in Figure \ref{fig:ap3b1ca} we get (note the paths $(a,c)$$\leftrightarrow$$x$$\leftrightarrow$$z$$\leftrightarrow$$b$) $\cind{ac}{b}{x}$ but $ac\not\ind b|zx$ (see supplement).  
On the other hand on the DAG in Figure \ref{fig:ap3b1c} clearly $\cind{ac}{b}{zx}$ but $ac\not\ind b|x$.  Similar examples for other four conditions can be drawn.

Theorem \ref{ap:3b} goes beyond the DAGs considered by \citet{sctsr3}. One example is considered in Figure \ref{fig:scts1}. Here $a\ind c$, $ac\ind z$ and both $\cind{ac}{b}{x}$ and $\cind{ac}{b}{xz}$ holds. 
Consequently, from Theorem \ref{ap:3b}, the relationship $\rhsq{a}{c}{b}\le\rhsq{a}{c}{bz^{\prime}}\le\rhsq{a}{c}{bz}$ follows.  Note that $z$ is not an \emph{ancestor} of $x$ but an \emph{ancestor} of $b$ and consequently, $zz^{\prime}\ind x$ also holds.  \citet{sctsr3} explicitly exclude conditioning vertices which are independent of $x$. 
  
\begin{cor}\label{cor:3b}
If $B=\emptyset$, Under all conditions of Theorem \ref{ap:3b} $(i)$, $\rhsq{a}{c}{z}=\rhsq{a}{c}{z^{\prime}}=\rho^2_{ac}=0$.  Under condition $(ii)$, $\rhsq{a}{c}{z}\ge\rhsq{a}{c}{z^{\prime}}\ge\rho^2_{ac}$.   
\end{cor}


\subsection{\bf Comparison between Theorems \ref{ap:1b} and \ref{ap:3b} for polytree models}
For polytree models, in view of Theorem \ref{ap:1b}, the conclusion of Theorem \ref{ap:3b} $(ii)$ is a bit counterintuitive. 
Note that, under $(ii)$, $\rhsq{a}{c}{x}=0$, which is same as in Theorem \ref{ap:1b}. However, unlike the latter, conditioning on vertices farther away produce a weaker squared correlation in this case.  
The difference seems to be that in Theorem \ref{ap:1b} $a\not\ind z$, but we assume $\cind{a}{z}{x}$.  In contrast, Theorem \ref{ap:3b} assumes that $a\ind z$, but in $(ii)$, the condition $a\ind z|x$ does not hold.
\begin{figure}[t]
\subfigure[\label{fig:scts1}]{
\input{figscts1.tex}
}\hspace{10pt}
\subfigure[\label{fig:cor4b}]{
\input{cor4b.tex}
}\hspace{10pt}
\subfigure[\label{fig:scts2}]{
\input{figscts2.tex}
}
\caption{\ref{fig:scts1} a DAG not considered by \citet{sctsr3}. From Theorem \ref{ap:3b}, it follows that $\rhsq{a}{c}{b}\le\rhsq{a}{c}{bz^{\prime}}\le\rhsq{a}{c}{bz}$.  \ref{fig:cor4b} A DAG to illustrate the contrast in the conclusion of Theorem \ref{ap:1b} and Theorem \ref{ap:3b} $(ii)$. Here $\rhsq{a}{c}{v}\ge\rhsq{a}{c}{u}\ge\rho^2_{ac}\ge\rhsq{a}{c}{y}\ge\rhsq{a}{c}{w}\ge\rhsq{a}{c}{x}=0$. 
From Theorem \ref{ap:1b}, on the DAG in \ref{fig:scts2} it follows that $\rhsq{a}{c}{b}\ge\rhsq{a}{c}{bz^{\prime}}\ge\rhsq{a}{c}{bz}$ always hold.}
\label{fig:extr}
\end{figure} 
As an illustration of this contrast we consider the graph in Figure \ref{fig:cor4b}. From Theorem \ref{ap:1b} and Corollary \ref{cor:3b} it follows that the relationship $\rhsq{a}{c}{v}\ge\rhsq{a}{c}{u}\ge\rho^2_{ac}\ge\rhsq{a}{c}{y}\ge\rhsq{a}{c}{w}\ge\rhsq{a}{c}{x}=0$ holds.  

Another such example can be constructed from the DAG in Figure \ref{fig:scts1}. We have argued above that from Theorem \ref{ap:3b} it follows that $\rhsq{a}{c}{b}\le\rhsq{a}{c}{bz^{\prime}}\le\rhsq{a}{c}{bz}$. In the DAG in Figure \ref{fig:scts2} the relation $a\ind c$ has been replaced by $\cind{a}{c}{x}$. 
From the rules of d-separation $\cind{a}{c}{bx}$, $\cind{ac}{z}{bx}$ and $\cind{acb}{z^{\prime}}{z}$ (see Definition \ref{defn:dconn}). Thus after conditioning on $b$, the Covariance matrix of $a$, $x$, $c$, $z$ and $z^{\prime}$ satisfies the conditions of Theorem \ref{ap:1b}. 
So the qualitative comparison holds, but in contrast to Figure \ref{fig:scts1}, it follows that $\rhsq{a}{c}{b}\ge\rhsq{a}{c}{bz^{\prime}}\ge\rhsq{a}{c}{bz}$.

\subsection{\bf Comparison between $\rhsq{a}{c}{x}$ and $\rhsq{a}{c}{Bz}$.}
If $z=x$, in Theorem \ref{ap:3b} in all case $a\ind Bcz$, so $\rhsq{a}{c}{z^{\prime}}=\rhsq{a}{c}{Bz}=\rhsq{a}{c}{B}=0$.  When $x\in V\setminus z$, comparison between $\rhsq{a}{c}{x}$ and $\rhsq{a}{c}{Bz}$ does not directly follow from Theorem \ref{ap:3b}.  
Under condition $(ii)$, $\cind{a}{c}{x}$, $0=\rhsq{a}{c}{x}\le\rhsq{a}{c}{Bz}$ for any $z$.  However, under the conditions $(i)$, $\rhsq{a}{c}{x}$ and $\rhsq{a}{c}{Bz}$ may not be qualitatively compared.  We show this fact in the following theorem.
\begin{theorem}\label{cl:3b2}
Suppose $a\ind z$, $c\ind az$, and $\cind{acz}{B}{x}$, then $\rhsq{a}{c}{Bz}\ge \rhsq{a}{c}{x}$, iff
\begin{equation*}\label{eq:firstcond}
\left(\sigma_{xx}+\frac{\sigma^2_{xz}}{\sigma_{zz}}\right)\Sigma_{xB}\Sigma^{-1}_{BB}\Sigma_{Bx}\ge \sigma^2_{xx},\text{  or equivalently  }\frac{\sigma_{xx}-\sigma_{xx|B}}{\sigma_{xx}}\ge\frac{\sigma_{zz|B}}{\sigma_{zz}}.
\end{equation*}
\end{theorem}

Theorems \ref{ap:1b}, \ref{ap:2b}, \ref{ap:3b} and \ref{cl:3b2} have a curious implication on polytree models.  Notice that in Theorems \ref{ap:1b} and \ref{ap:2b} the vertex $z$ is in the set of \emph{descendants} of vertex $x$ (see Figures \ref{fig:figap1b2} and \ref{fig:ap2b}), 
whereas in Theorem \ref{ap:3b}, $z$ may be a \emph{parent} of $x$. 
The curious fact is that, on a polytree the squared 
partial correlations given the descendants of $x$ cannot be compared with the squared partial correlations given the parents (or more generally given the \emph{ancestors of the parents} of $x$). 
Furthermore, the behaviour of $\rhsq{a}{c}{x}$ is a continuation of the behaviour of squared partial correlations given its descendants.  
In other words, on polytrees, conditioning on the vertices ``above'' the path has different nature than conditioning on the vertices ``below'' or ``on'' the path.     

\begin{figure}[t]
\subfigure[\label{fig:compgraph}]{
\input{compgraph.tex}
}\hspace{10pt}
\subfigure[\label{plot:compgraph}]{
\includegraphics[height=2in,height=2in]{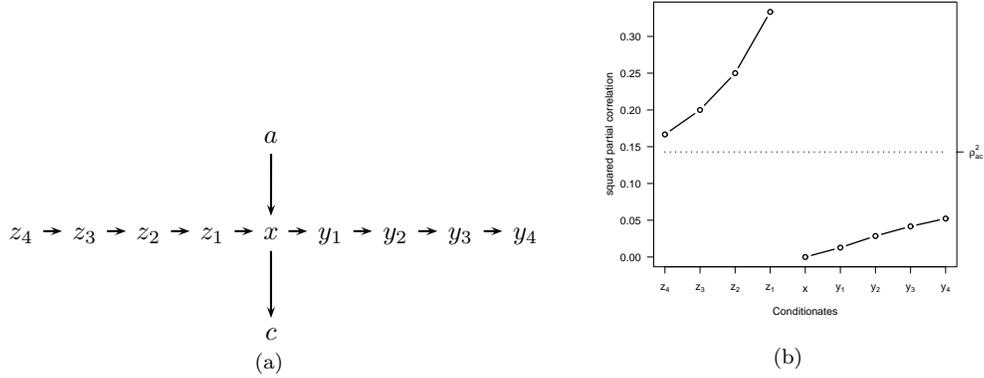}
}
\caption{\ref{fig:compgraph} A polytree and \ref{plot:compgraph} the value of $\rhsq{a}{c}{i}$ for $i\in\{\emptyset,z_4,z_3,z_2,z_1,x,y_1,y_2,y_3,y_4\}$. Each parameter is fixed at $1$. \ref{plot:compgraph} illustrates the discontinuous drop in $\rhsq{a}{c}{i}$ as we move from $z_1$ to $x$ along the $z_4$ to $y_4$ path. }
\label{fig:comps}
\end{figure}

We present an illustrative example in Figure \ref{fig:comps}.  We consider the polytree in Figure \ref{fig:compgraph}. In Figure \ref{plot:compgraph} we plot the values of $\rhsq{a}{c}{i}$ for $i\in\{\emptyset,z_4,z_3,z_2,z_1,x,y_1,y_2,y_3,y_4\}$. All parameter values are fixed at $1$. As predicted from Theorem \ref{ap:3b} the squared partial correlation increases from $i=z_4$ to $i=z_1$ and from Corollary \ref{cor:3b} each of them are larger than $\rho^2_{ac}$. 
However, From Theorem \ref{ap:2b}, $\rhsq{a}{c}{i}$ increases as we move from $x$ to $y_4$ and each of them are smaller that $\rho^2_{ac}$. Thus the squared partial correlation drops discontinuously as we move from $z_1$ to $x$ along the $z_4$ to $y_4$ path.     
 
\subsection{\bf Further generalisations on comparison with fixed correlates}
Suppose $Z_1=\{z_{11},z_{12},\ldots,z_{1n}\}$ and $Z_2=\{z_{21},z_{22},\ldots,z_{2n}\}$ are two conditionates of cardinality $n$. Then for fixed correlates $a$ and $c$, one can write:
\begin{equation}\label{eq:factmain}
\frac{\rhsq{a}{c}{Z_1}}{\rhsq{a}{c}{Z_2}}=\prod^n_{i=1}\frac{\rhsq{a}{c}{z_{21},z_{22},\ldots,z_{2(i-1)},z_{1i},z_{1(i+1)},\ldots,z_{1n}}}{\rhsq{a}{c}{z_{21},z_{22},\ldots,z_{2(i-1)},z_{2i},z_{1(i+1)},\ldots,z_{1n}}}
\end{equation}  
Clearly $\rhsq{a}{c}{Z_1}\le\rhsq{a}{c}{Z_2}$ holds if each factor in the R.H.S. of \eqref{eq:factmain} is bounded by $1$.  

Note that in each factor in \eqref{eq:factmain} the conditionate in the numerator and the denominator differ only in one element.  
Thus in order to qualitatively compare $\rhsq{a}{c}{Z_1}$ and $\rhsq{a}{c}{Z_2}$it is sufficient to find a $x_i$ for each factor such that $z_{1i}$ and $z_{2i}$ satisfy the conditions of one of the Theorems $2$ - $4$, possibly with $B\subseteq\{z_{21},z_{22},\ldots,z_{2(i-1)},z_{1(i+1)},\ldots,z_{1n}\}$ whenever necessary. 

Using the factorisation in \eqref{eq:factmain} and Theorems $2$ - $4$, structural and path based rules for comparison may be postulated for several graphical models.  The choice of $x_i$ and these path based rules depend on the structure of association of the whole vector $V$.  We consider the tree models below.

\section{Application to tree models}\label{sec:tree}
Let $\gr{G}{V}{E}$ be a tree with vertex set $V$ and edge set $E$.  For vertices $x\in V$ and $y\in V$, \path{x}{y} denote the unique path joining $x$ and $y$, which we define as:
\begin{equation}
\begin{aligned}
\epath{x}{y}=\{x=v_1,v_2,\ldots,v_{k-1},v_k=y\text{ such that there is an edge between}&\\\nonumber
\text{$v_i$ and $v_{i+1}$, for each $i=1,2,\ldots,k-1$}\}.&
\end{aligned} 
\end{equation}  

Notice that, by the above definition \path{x}{y} is a subset of $V$ which contains the end points $x$ and $y$. Since $G$ is a tree, it has only one connected component and therefore any two vertices $x$ and $y$ are connected by an unique \path{x}{y}. 

\begin{defn}\label{defn:sep}
 Two vertices $a$ and $c$ on an undirected graph $G$ is said to be \emph{separated} given a subset $Z$ of $V\setminus\{a,c\}$ if each path \path{}{} between $a$ and $c$ intersects $Z$. 
Two subsets $A$ and $C$ of $V$ are separated given $Z\subseteq V\setminus(A\cup C)$ if $Z$ separates each $a\in A$ from each $c\in C$.    
Two subset $A$ and $C$ of $V$ are connected given a subset $Z$ if they are not separated given $Z$.
\end{defn}

Clearly on a tree $a$ and $c$ are separated given each $x\in\epath{a}{c}\setminus\{a,c\}$.  On the other hand since any two vertices $a$ and $c$ are connected by an unique path, $a$ and $c$ cannot be separated given the $\emptyset$. 

The \emph{separation criterion} described above associates a set of conditional independence relations with $G$.  This set is described by a collection of \emph{triples}.
\begin{equation}
\indrels{G}=\left\{\trip{T_1}{T_2}{T_3},\text{where $T_1\dot{\cup}T_2\dot{\cup}T_3\subseteq V$ such that $\cind{T_1}{T_2}{T_3}$}\right\}. 
\end{equation}
The association of the separation criterion with $\indrels{G}$ can be described as follows:
\[\trip{T_1}{T_2}{T_3}\Leftrightarrow \text{ $T_1$ is separated from $T_2$ given $T_3$ in $G$.}
\] 
If $V\sim N\left(0,\Sigma\right)$, then $\Sigma$ satisfies all conditional independence relationships in $\indrels{G}$.  This implies that if $\Lambda=\Sigma^{-1}$, for each $\trip{T_1}{T_2}{T_3}\in\indrels{G}$, $\Lambda_{T_1T_2}=0$.

We now define formal operation of conditioning for independence model $\indrels{G}$, on subsets of $V$.

\begin{defn}\label{defn:condrels}
An independence model $\indrels{G}$ {\it after conditioning on a subset} $Z$ is the set of triples defined as follows:
\begin{equation}
\cm{\indrels{G}}{Z}{\empty}\; \equiv\; \Bigl\{ \trip{T_1}{T_2}{T_3} 
\; \Bigr|\; \Bigl.
\trip{T_1}{T_2}{T_3\cup Z}
\in {\indrels{G}};\;\, (T_1\cup T_2\cup T_3)\cap Z = \emptyset \Bigr\}.
\end{equation}
\end{defn}
Thus if $\indrels{G}$ contains the independence relations satisfied by a $N\left(0,\Sigma\right)$ on $G$, then
$\cm{\indrels{G}}{Z}{\empty}$ constitutes the subset of independencies holding among the variables in $Z^c=V\setminus Z$, after conditioning on $Z$.  Let $G_{Z^c}$ be the subgraph of $G$ with vertex set $Z^c$ and edge set consisting of all edges in $E$ between the vertices in $Z^c$.  The following Lemma makes the connection between $\cm{\indrels{G}}{Z}{\empty}$ and $\mathfrak{I}\left(G_{Z^c}\right)$.

\begin{lem}\label{lem:ugcond}
Suppose $\gr{G}{V}{E}$ is a tree. Let $a,c$ be two distinct vertices, $Z\subseteq V\setminus\{a,c\}$ and $Z^c=V\setminus Z$.  Then  
\begin{equation}
\cm{\indrels{G}}{Z}{}=\mathfrak{I}\left(G_{Z^c}\right).
\end{equation} 
\end{lem}

Lemma \ref{lem:ugcond} holds for any UG. It implies that the conditioning on $Z$ does not add or delete any edge in $G_{Z^c}$, so if $G$ is tree $\cm{\indrels{G}}{Z}{}$ can be represented by a forest.  The inverse of conditional covariance matrix of $Z^c$ given $Z$ is simply $\Lambda_{Z^cZ^c}$.

Separation ensures conditional independence, but if even if the separation fails the corresponding conditional covariance can still be zero  (implying conditional independence for Gaussian random variables) because of the parameter values. However, Theorem \ref{ap:1b} is still valid in these cases. 

For a fixed conditionate the rules for comparing squared partial correlations on trees follows easily from Theorem \ref{lem:maincomp1} and the separation criterion.
\begin{theorem}\label{thm:cond}
Suppose that, on a Gaussian tree $G$, the vertices $a$, $c$, $c^{\prime}$ are such that $c\in\epath{a}{c^{\prime}}$. Then for any $Z\subseteq V$, $\rhsq{a}{c^{\prime}}{Z}\le\rhsq{a}{c}{Z}$. 
\end{theorem} 

For fixed correlates $a$ and $c$ and two sets $Z_1$ and $Z_2$ of cardinality more than one, $\rhsq{a}{c}{Z_1}$ and $\rhsq{a}{c}{Z_2}$ can be compared qualitatively. The following result describes a sufficient condition.

\begin{theorem}\label{thm:ugcorr}
Let $\gr{G}{V}{E}$ be a Gaussian tree. Suppose $a$ and $c$ are two vertices on $G$ and $Z_1$ and $Z_2$ are two subsets of $V$ such that $\cind{ac}{Z_2}{Z_1}$.  Then $\rhsq{a}{c}{Z_1}\le\rhsq{a}{c}{Z_2}$.  
\end{theorem}

From the separation criterion described above, it follows that the vertices $a$ and $c$ separated from $Z_2$ given $Z_1$ implies $\cind{ac}{Z_2}{Z_1}$ and therefore $\rhsq{a}{c}{Z_1}\le\rhsq{a}{c}{Z_2}$.  
The following Corollary gives the corresponding sufficient condition in terms of paths: 

\begin{cor}\label{cor:ugcorr} 
Suppose $Z_1$ and $Z_2$ are two subsets of $V$, such that for each vertex $z_2\in Z_2$, the both paths \path{a}{z_2} and \path{c}{z_2} intersect $Z_1$, then $\rhsq{a}{c}{Z_1}\le \rhsq{a}{c}{Z_2}$. 
\end{cor}

Notice that, Theorem \ref{thm:ugcorr} is more general than Corollary \ref{cor:ugcorr}, the Theorem covers the cases when the conditional independence holds due to the choices of parameters as well.  The result in Theorem \ref{thm:ugcorr} is also complete in the following sense.  
\begin{theorem}\label{thm:cmplt}
Suppose $\gr{G}{V}{E}$ is a Gaussian tree. Let $Z_1,Z_2\subseteq V$ such that $ac\not\ind Z_2|Z_1$ and $ac\not\ind Z_1|Z_2$.  Further, suppose that $\left(Z_1\cup Z_2\right)\cap\epath{a}{c}=\emptyset$.  
Then there exists $\Sigma_1$ such that $\rhsq{a}{c}{Z_1}>\rhsq{a}{c}{Z_2}$ and $\Sigma_2$ such that $\rhsq{a}{c}{Z_2}>\rhsq{a}{c}{Z_1}$. 
\end{theorem}

Finally, Theorem \ref{thm:cond} and the Corollary \ref{cor:ugcorr} can be combined to a general rule for comparing squared partial correlation on trees.

\begin{cor}\label{cor:final}
Suppose $a$, $c$, $c^{\prime}$ are three vertices on a Gaussian tree $G$ and $Z$, $Z^{\prime}$ are two subsets of the vertex set $V$. Further, assume that $c\in\epath{a}{c^{\prime}}$ and the vertices $a$ and $c^{\prime}$ are separated from $Z$ given $Z^{\prime}$.  Then $\rhsq{a}{c^{\prime}}{Z^{\prime}}\le \rhsq{a}{c}{Z}$.
\end{cor} 
 

\section{Application to polytree models and model selection}
\begin{figure}[t]
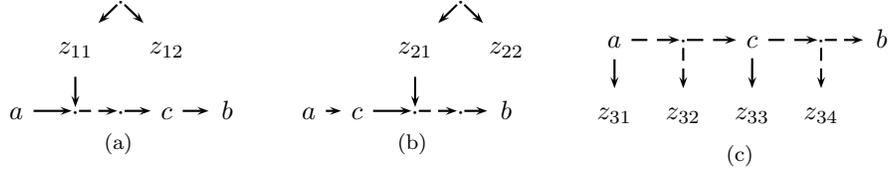

\begin{center}
\subfigure[\label{fig:modsel1}]{
\input{dag_sel1.tex}}\hspace{.1in}
\subfigure[\label{fig:modsel2}]{
\input{dag_sel2.tex}
}\hspace{.1in}
\subfigure[\label{fig:modsel3}]{
\input{dag_sel3.tex}
}
\end{center}
\caption{Examples of polytrees satisfying the conditions of Theorem \ref{thm:modsel} below.  In each, $a\in an(c)$ and $c\in an(b)$.  In \ref{fig:modsel1} $z_{11}$ and $z_{12}$ satisfy condition $1.$ and $\rhsq{a}{c}{b}<\rhsq{a}{c}{bz_{12}}<\rhsq{a}{c}{bz_{11}}$ (from Theorem \ref{ap:3b}. (ii)).  
In \ref{fig:modsel2} $z_{21}$ and $z_{22}$ satisfy condition $2.$. So $\rhsq{a}{c}{bz_{21}}<\rhsq{a}{c}{bz_{22}}<\rhsq{a}{c}{b}$ (see Theorem \ref{ap:1b} and Figure \ref{fig:scts2}). Each $z_{3k}$, $k=1,\ldots,4$, in \ref{fig:modsel3} satisfy condition $2.$, 
ie. $\rhsq{a}{c}{z_{3k}}<\rhsq{a}{c}{b}$. Note that, $b$ cannot be in $an(z)$, otherwise $\cind{ac}{z}{b}$ and $\rhsq{a}{c}{b}= \rhsq{a}{c}{bz}$.}
\label{fig:modsel}
\end{figure}

A polytree is a DAG such that if we substitute all its directed edges with undirected ones, the resulting graph (ie. its skeleton) would be a tree.  Thus on a polytree two vertices $x$ and $y$ can have at most one path \path{x}{y} connecting them.  
Here, on a connecting path we disregard the direction of the individual edges.

A vertex $y$ is an ancestor of a vertex $x$, if either $y=x$ or $x$ can be reached from $y$ by following the arrowheads of a directed path (ie. the path $y\rightarrow v_1\rightarrow v_2 \dashrightarrow v_k\rightarrow x$ exits).  The collection of all ancestors of $x$ is denoted by $an(x)$.   Furthermore, for a set of vertices $X$ we define $an(X)=\cup_{x\in X}an(x)$.

\begin{theorem}\label{thm:modsel} Suppose that on a Gaussian polytree $a\ne c\ne b$, $a\in an(c)$ and $c\in an(b)$.  Further let, for some vertex $z$, $\rhsq{a}{c}{bz}\ne \rhsq{a}{c}{b}$.  Then\begin{enumerate}
\item $\rhsq{a}{c}{bz}>\rhsq{a}{c}{b}$, iff $a\ind z$ and $c\not\ind z$.
\item $\rhsq{a}{c}{bz}<\rhsq{a}{c}{b}$ iff either $c\ind z$ or $a\not\ind z$.
\end{enumerate} 
\end{theorem} 

The condition $\rhsq{a}{c}{bz}\ne \rhsq{a}{c}{b}$ is required in Theorem \ref{thm:modsel}. This implies $ac\not\ind z\mid b$. So $b\not\in an(z)$.  It can further be shown (see the proof) that the polytree structure implies $ac\ind z$ iff $c\ind z$.  
Thus the right hand side of Condition $2.$ above equivalently means that either both $a$ and $c$ are independent of $z$ or none of them are independent of $z$. Examples of graphs satisfying the conditions $1.$ and $2.$ can be found in Figure \ref{fig:modsel}.


\begin{figure}[t]
\begin{center}
\includegraphics[width=9cm,keepaspectratio]{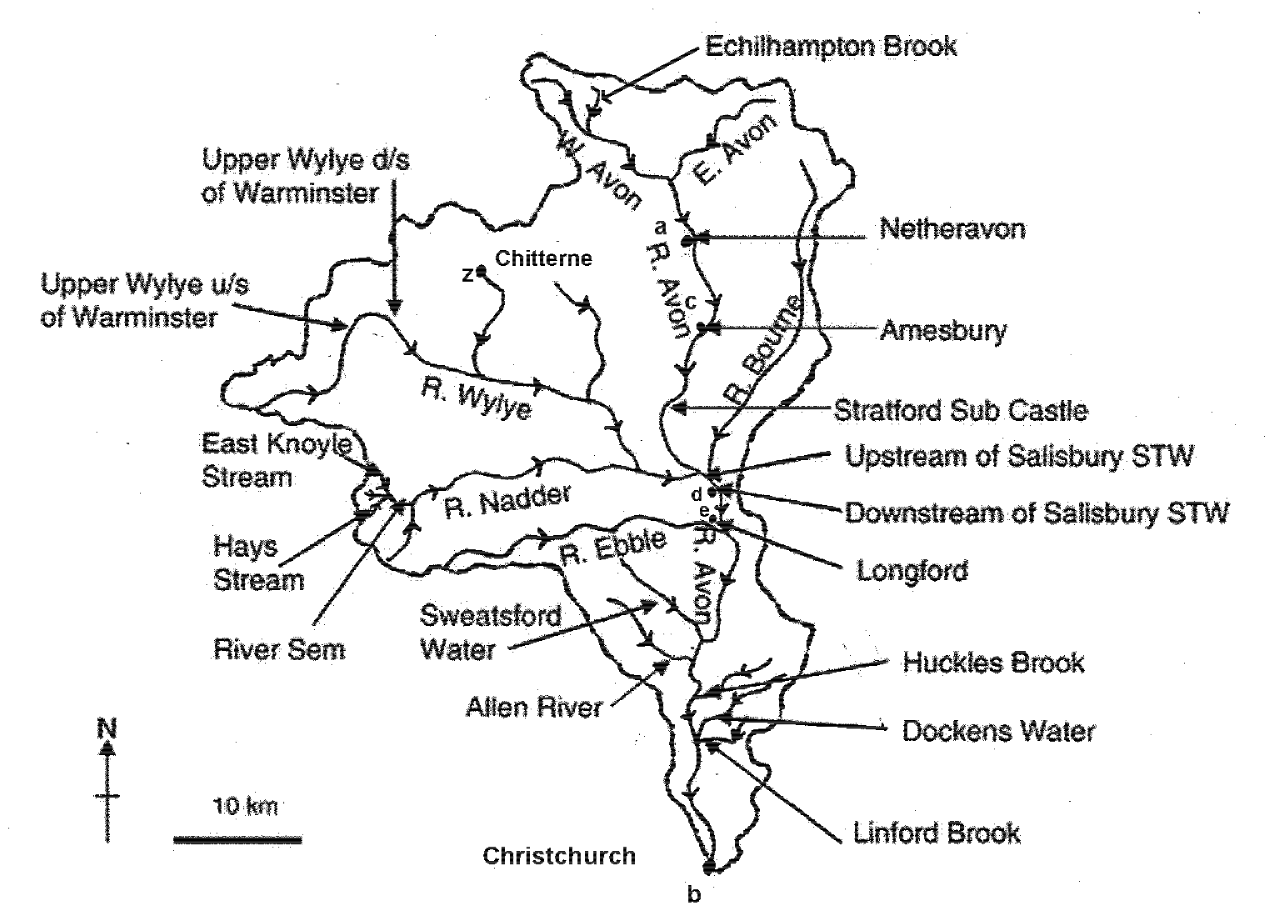}
\end{center}
\caption{An illustration of the results in Theorem \ref{thm:modsel} on the river network of Avon river, Hampshire, England (obtained from \citet{jarvie:etal:2005}).}
\label{fig:avon}
\end{figure}

Theorem \ref{thm:modsel} has applications in model selection. An example occurs in the mapping of river flow networks.  
Figure \ref{fig:avon} \citep{jarvie:etal:2005} presents a schematic diagram of the network of the Avon basin in Hampshire, England. Suppose that it is known that none of the rivers involved have a distributary.  Clearly the network, with the direction of the water flow form a polytree.  
Measurements can be taken at points $a$ (Netheravon), $b$ (Christchurch), $c$ (Amesbury), $d$ (Downstream of Salisbury STW), $e$ (Longford) and $z$ (Chitterne). 
However, because of practical considerations we suppose that the measurements are taken when the water level at Christchurch ($b$) touches certain levels. Lets assume $\rhsq{a}{x}{b}\ne\rhsq{a}{x}{bz}$ for $x=c,d,e$.  We want to know where does the stream from $z$, ie. Chitterne meets river Avon.  

It is clear that since the observations are all conditional on the water level at $b$, in the data neither $z\not\ind a$ nor $z\not\ind c$.  However, from Theorem \ref{ap:1b}, see also Figure \ref{fig:scts2} and Theorem \ref{ap:3b} it follows that $\rhsq{a}{c}{bz}<\rhsq{a}{c}{b}$, $\rhsq{a}{d}{bz}>\rhsq{a}{d}{b}$ and $\rhsq{a}{e}{bz}>\rhsq{a}{e}{b}$.  
From Condition $2.$ of Theorem \ref{thm:modsel} it follows that either both $a$ and $c$ are independent of $z$ or none of them are. On the other hand, Condition $1.$ implies that $a\ind z$ but $d$ and $e$ are not independent of $z$.  
If none of $a$ and $c$ are independent of $z$, the point $z$ must be on a distributary stream or on a tributary which meets Avon north of $a$ (Netheravon).  However, by assumption there is no distributary stream.  Furthermore, if the tributary from $z$ meets Avon somewhere north of $a$, by Theorem \ref{ap:1b} both $\rhsq{a}{d}{bz}<\rhsq{a}{d}{b}$ and $\rhsq{a}{e}{bz}<\rhsq{a}{e}{b}$ must hold. 
This is a contradiction. Thus $ac\ind z$ must hold.  So from Theorem \ref{thm:modsel} we see that the stream from Chitterne ie. $z$ meets Avon somewhere between Amesbury ie. $c$ and Downstream of Salisbury STW ie.$d$. 


\section{Necessity of the conditional independence relationships}
In the above sections we postulated some sufficient conditional independence relationships under which some squared conditional correlations can be qualitatively compared. It is not known if these relationships are necessary as well.  
It is possible that qualitative comparison would hold under different sets of conditions. However the conditions in any set of relationships cannot be reduced. In this section we show this fact using various counterexamples.

In each counter-example, unless otherwise stated, set all parameters ie. the regression coefficients and the node specific conditional variances are set to $1$.
\subsection{Comparison with a fixed conditionate}
\begin{figure}[t]
\parbox{.4\columnwidth}{
\begin{picture}(150,200)(0,0)
\thicklines
\put(130,170){\vector(-1,-1){55}}
\put(20,170){\vector(1,-1){55}}
\put(75,115){\vector(0,-1){65}}
\put(132.5,172.5){$a$}
\put(72.5,122.5){$c$}
\put(15,172.5){$x$}
\put(72.5,40){$z$}
\end{picture}
}\ \hspace{.5in} \
\parbox{.4\columnwidth}{
\resizebox{.4\columnwidth}{.36\columnwidth}{\includegraphics{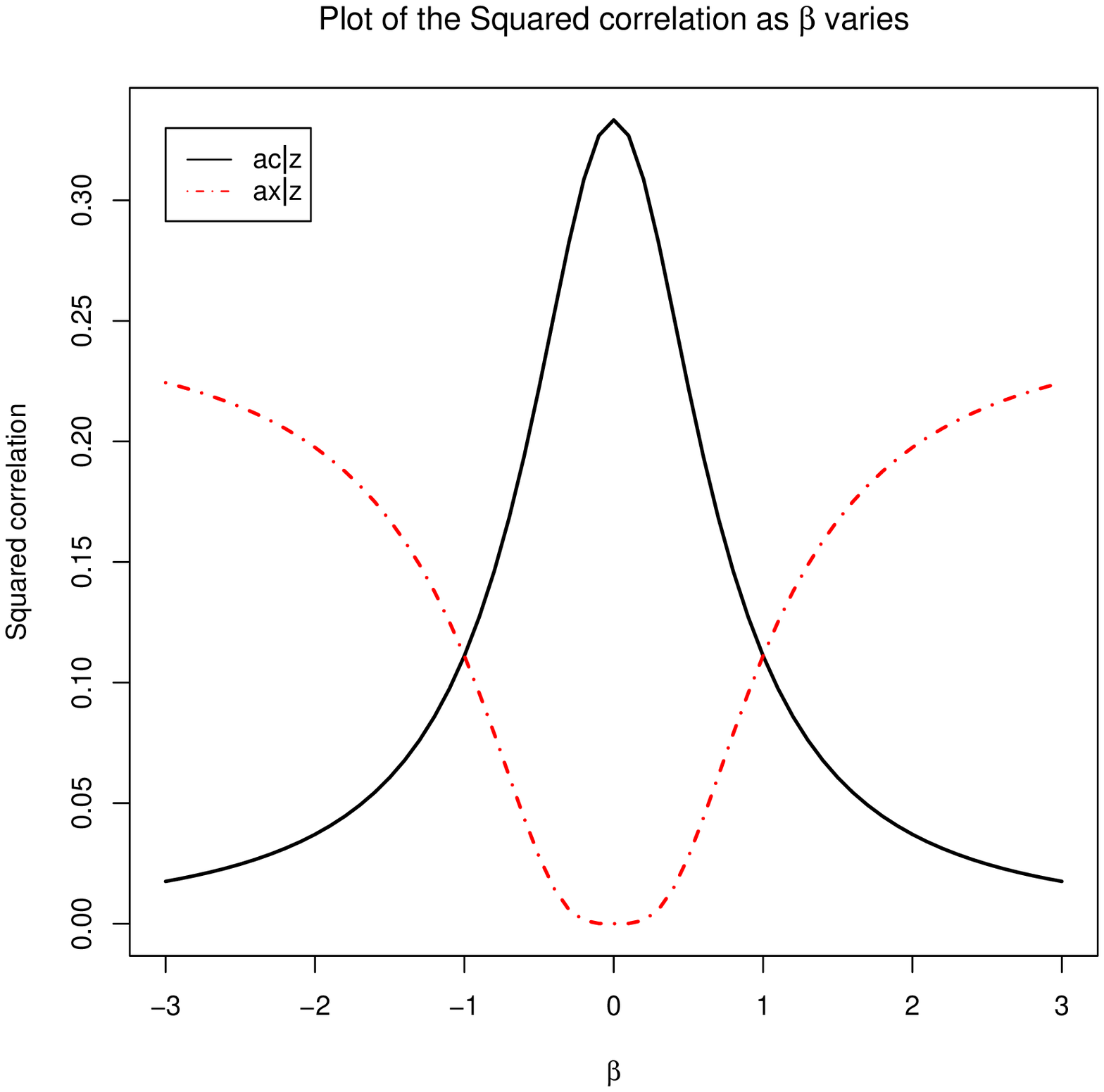}}
\resizebox{.4\columnwidth}{.36\columnwidth}{\includegraphics{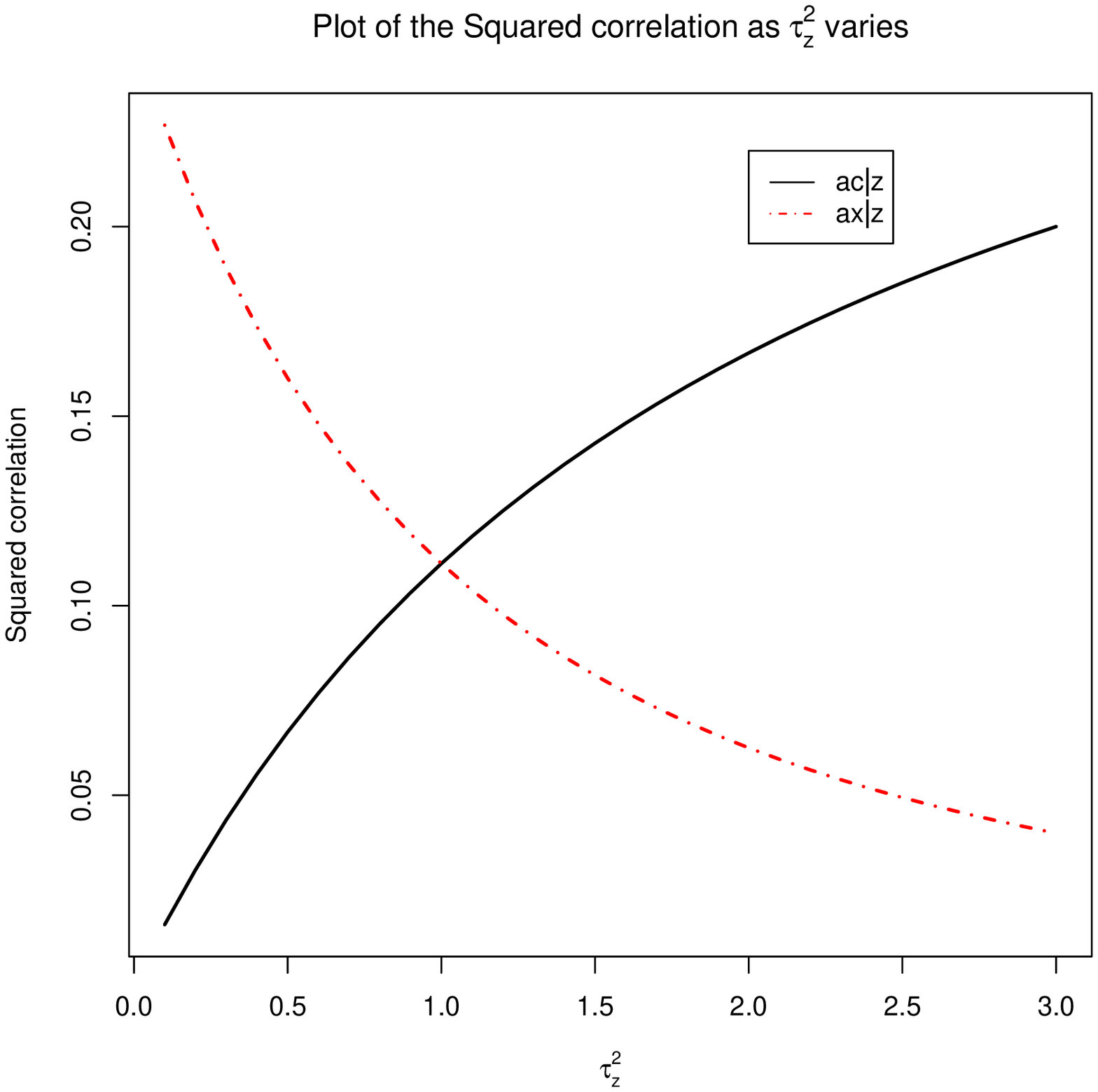}}
}
\caption{Plot of $\rhsq{a}{c}{z}$ and $\rhsq{a}{x}{z}$ for the graph on the
left with $\beta_{zc}$ and $\tau^2_z$.}
\label{figcomp1:2}
\end{figure}

We consider the graph in Figure \ref{figcomp1:2}.  Note that, $c$ is a collider on the \path{a}{x} and $z$ is a child of $c$. Thus, from the laws of d-separation $x$ is not d-separated from $a$ given $c$ and $z$. 
Under our choice of parametrisation clearly $a\not\ind x\mid cz$.  In the plots to the right of Figure \ref{figcomp1:2} we change respectively $\beta_{cz}$ and $\tau^2_z$ and keep other parameters fixed. It is clear from the plots that $\rhsq{a}{c}{z}$ and $\rhsq{x}{c}{z}$ cannot be qualitatively compared. This shows the condition of Theorem \ref{lem:maincomp1} cannot be relaxed.    
    
\subsection{Comparison with fixed correlates}
We only consider the necessity of the conditions of Theorems \ref{ap:1b} and \ref{ap:2b} here. The examples for Theorem \ref{ap:3b} are similar.

\begin{figure}[t]
\hspace{1in}\subfigure[\label{fig:c1g}]{
$\psmatrix[colsep=.6in,rowsep=.6in]
&x&\\
a& &c\\
&z_1&\\
&z_2&
\endpsmatrix
$
\psset{nodesep=3pt}
\ncline{->}{1,2}{2,1}
\ncline{->}{1,2}{2,3}
\ncline{->}{1,2}{3,2}\lput*{0}{}
\ncline[linestyle=dashed]{->}{2,3}{3,2}\Aput{\small{$\beta_{z_1c}$}}
\ncline{->}{3,2}{4,2}
\ncline[linestyle=dashed]{->}{2,3}{2,1}\lput{:0}{\rput{N}(.6,.4){$\beta_{ac}$}}
\ncline[linestyle=dashed]{->}{2,1}{4,2}\Bput{$\beta_{z_2a}$}
}\hfill
\subfigure[\label{fig:c1e1}]{
\resizebox{3in}{3in}{\includegraphics{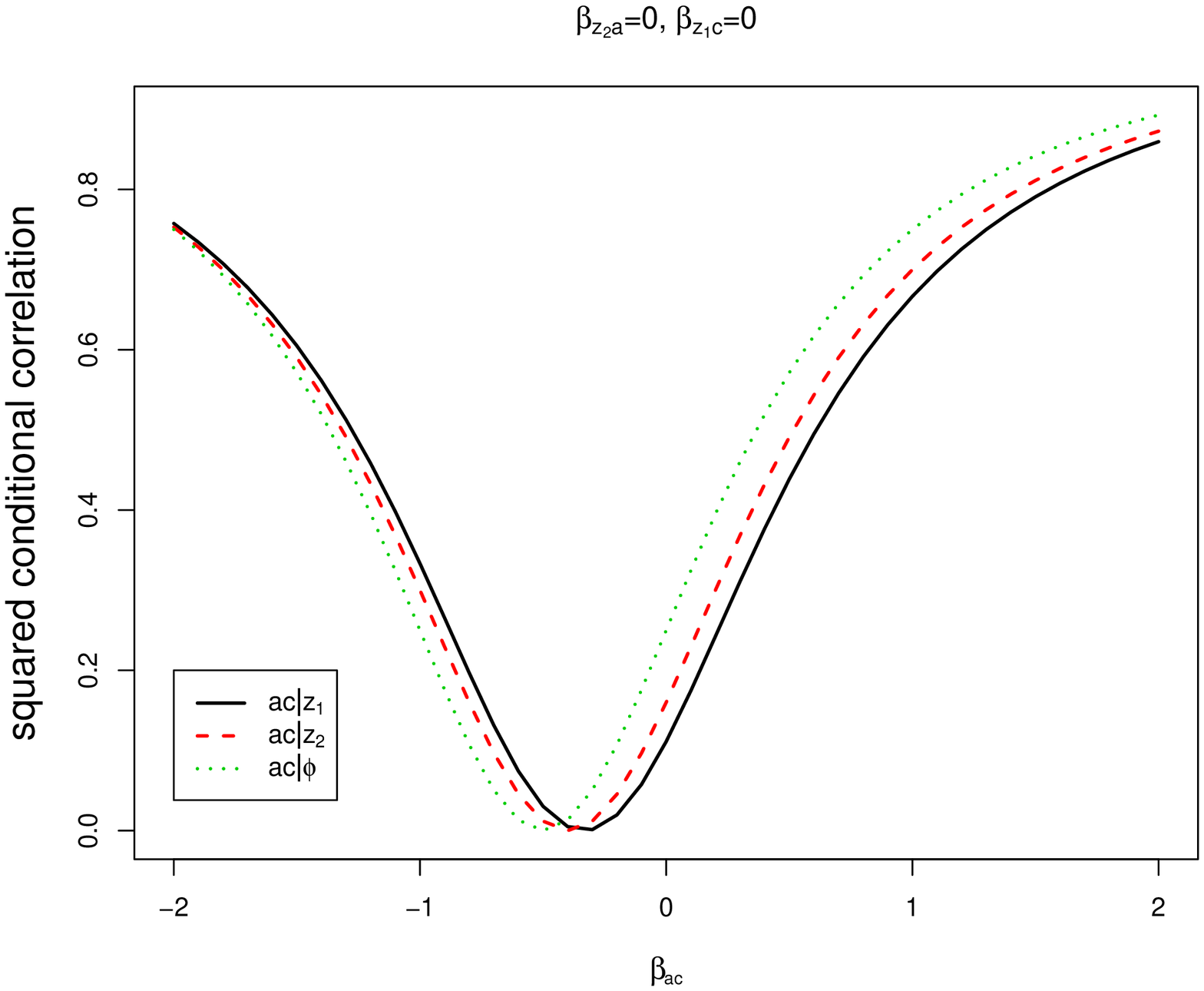}}
}\hfill
\subfigure[\label{fig:c1e2}]{
\resizebox{3in}{3.1in}{\includegraphics{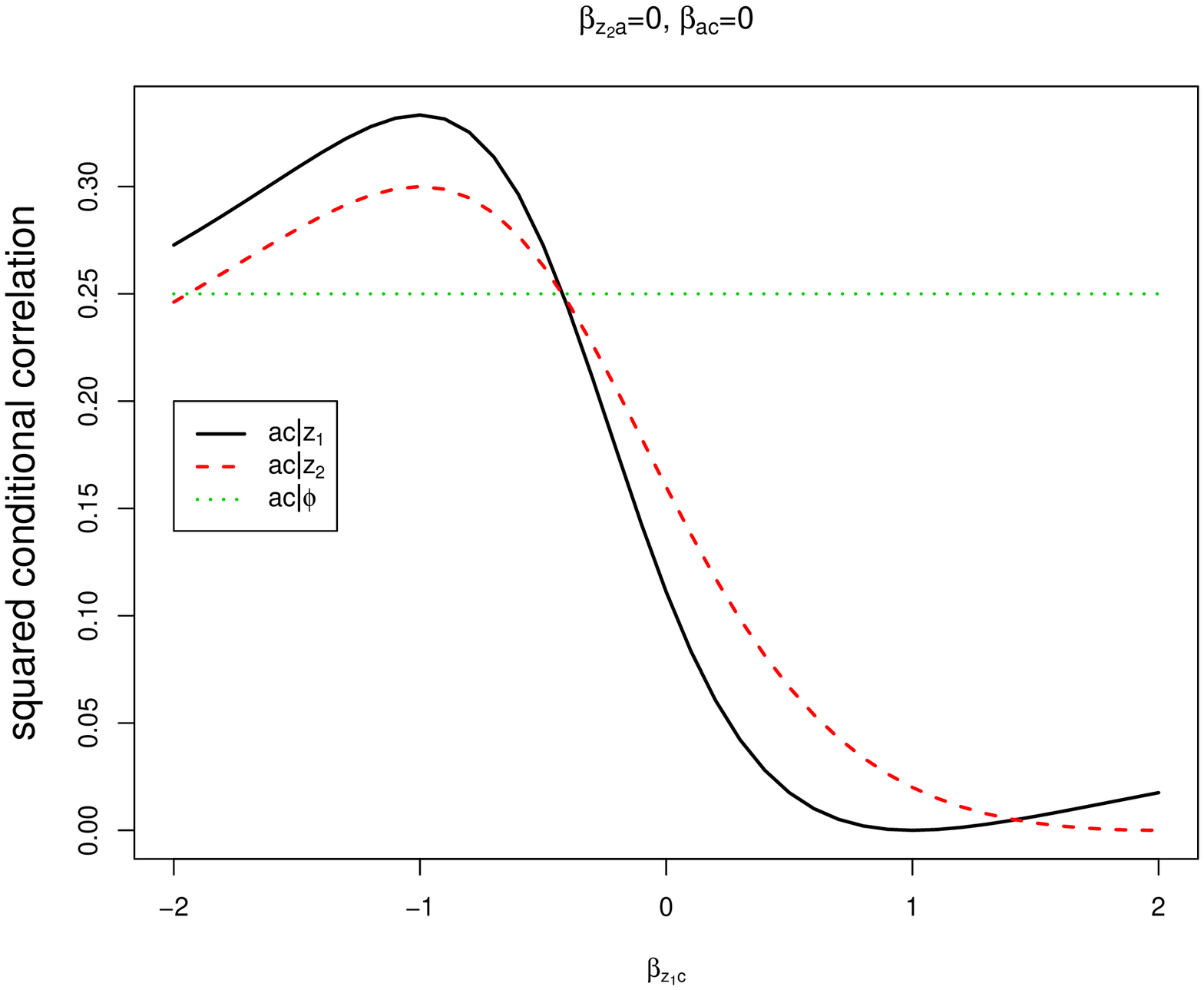}}
}\hfill
\subfigure[\label{fig:c1e3}]{
\resizebox{3in}{3.1in}{\includegraphics{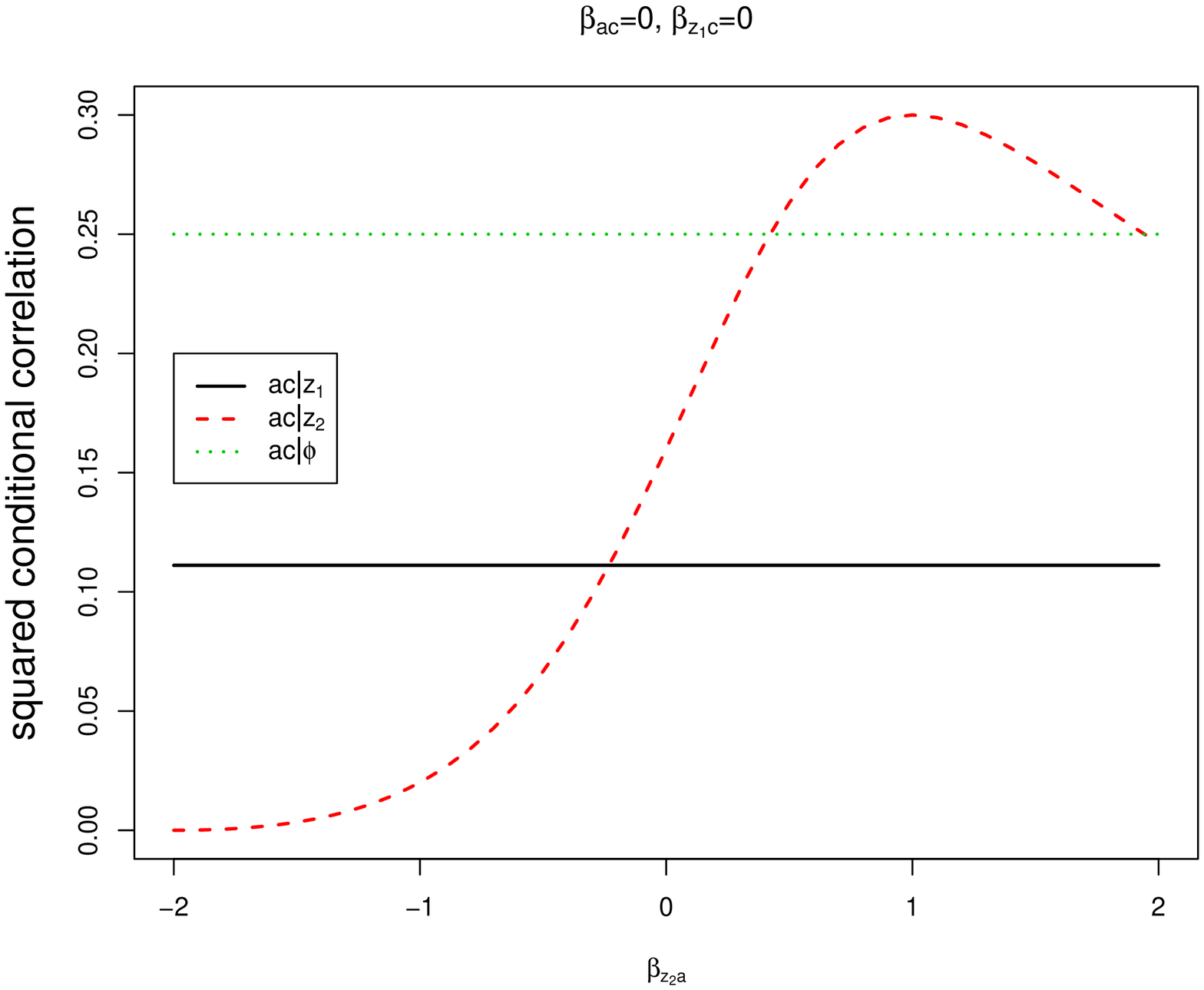}}
}
\caption{The directed acyclic graph described in Example 1. Clearly there is no ordering between $\rho^2_{ac}$, $\rhsq{a}{c}{z_1}$ and $\rhsq{a}{c}{z_2}$.}
\label{fig:1}
\end{figure}

\begin{figure}[t]
\hspace{1in}\subfigure[\label{fig:c2g}]{
$\psmatrix[colsep=.6in,rowsep=.6in]
a&&c\\
&x&\\
b&z_1&\\
&z_2&
\endpsmatrix
$
\psset{nodesep=3pt}
\ncline{->}{1,1}{2,2}
\ncline{->}{1,3}{2,2}
\ncline[linestyle=dashed]{->}{1,1}{1,3}\Aput{$\beta_{ca}$}
\ncline[linestyle=dashed]{->}{1,3}{3,2}\Aput{\small{$\beta_{z_1c}$}}
\ncline{->}{2,2}{3,2}
\ncline{->}{3,2}{4,2}
\ncline{->}{2,2}{3,1}
\nccurve[ncurv=.4,angleB=45,angleA=315,linestyle=dashed]{->}{1,3}{4,2}\Aput{$\beta_{z_2c}$}
}\hfill
\subfigure[\label{fig:c2e1}]{
\resizebox{3in}{3in}{\includegraphics{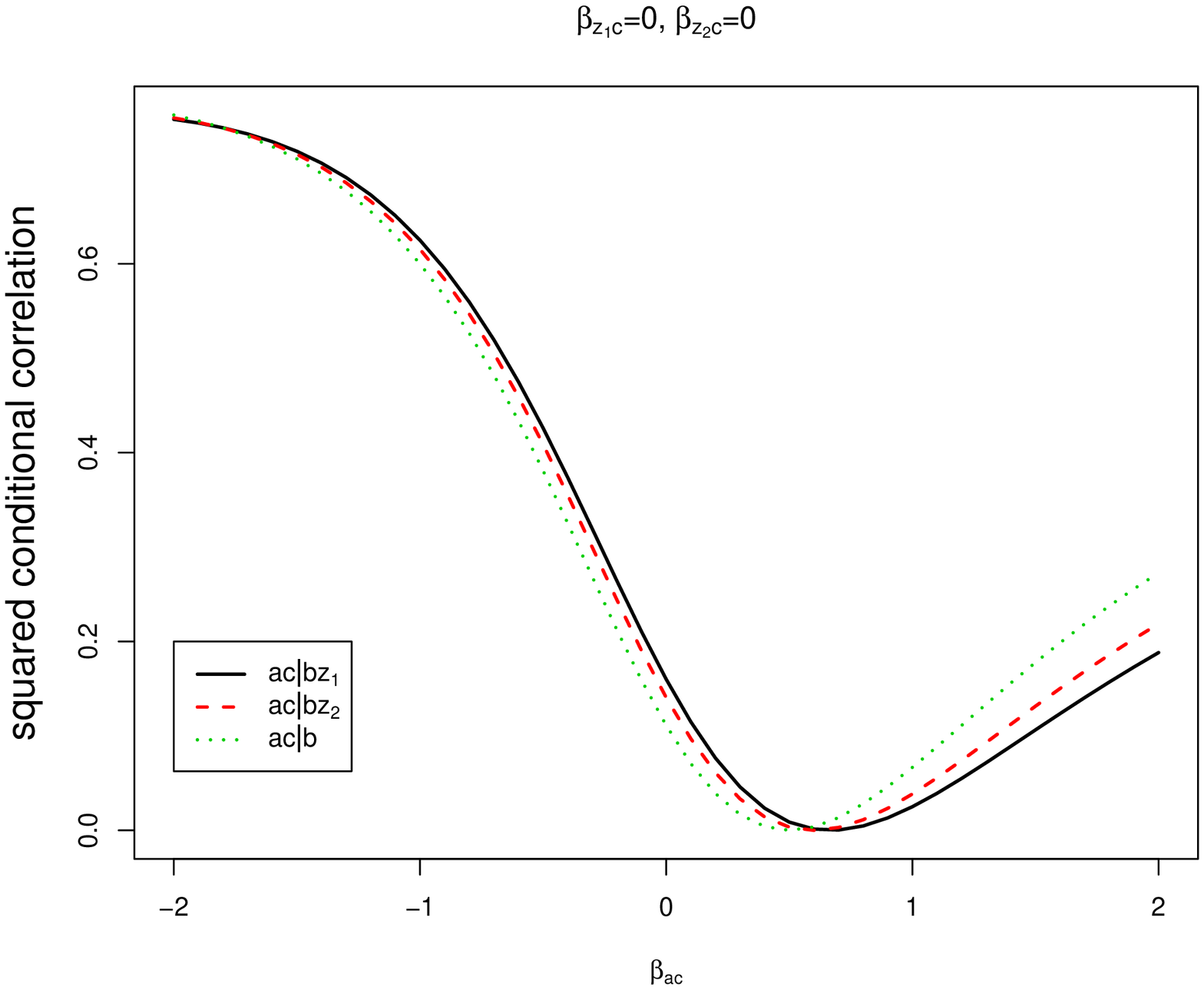}}
}\hfill
\subfigure[\label{fig:c2e2}]{
\resizebox{3in}{3in}{\includegraphics{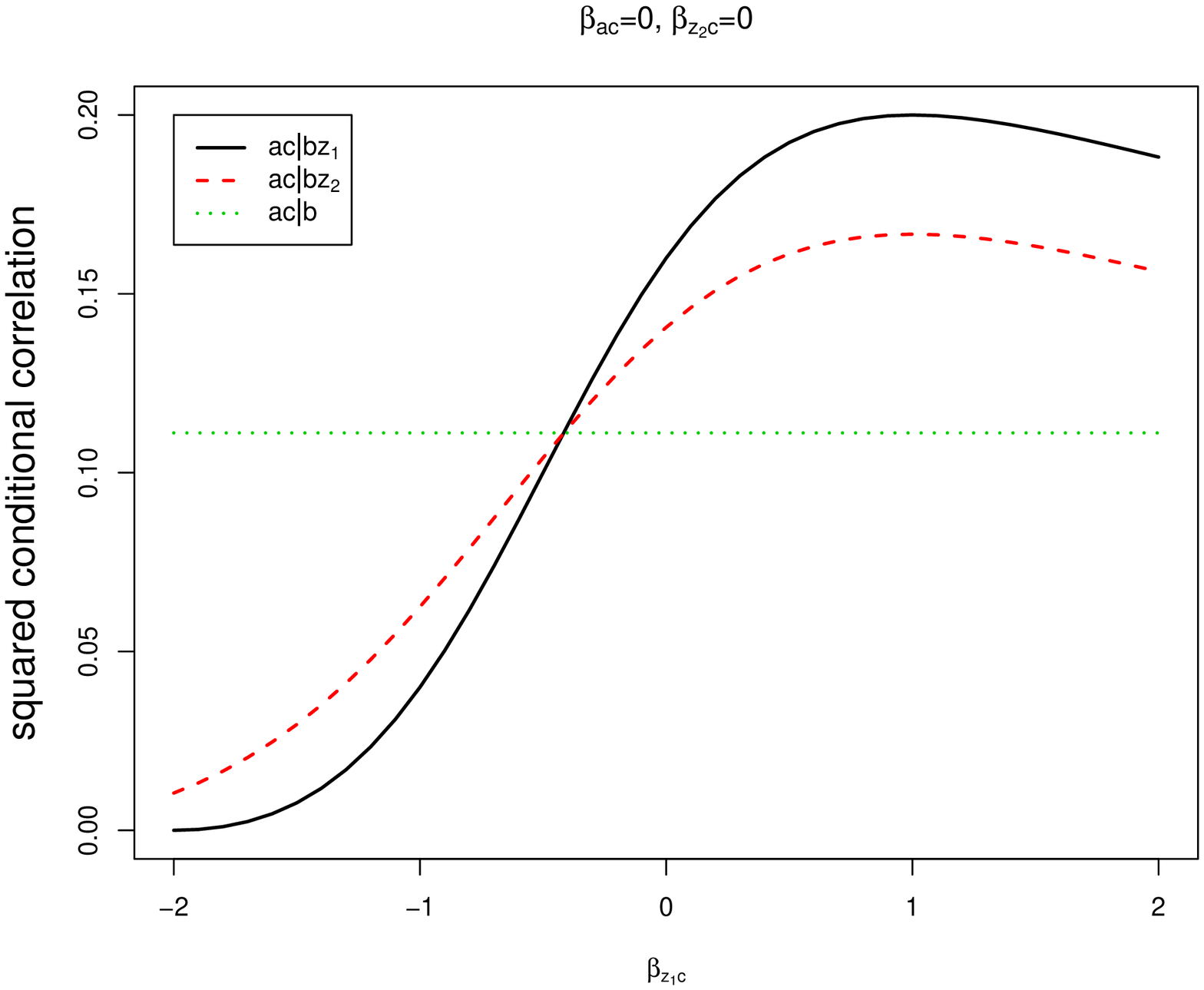}}
}\hfill
\subfigure[\label{fig:c2e3}]{
\resizebox{3in}{3in}{\includegraphics{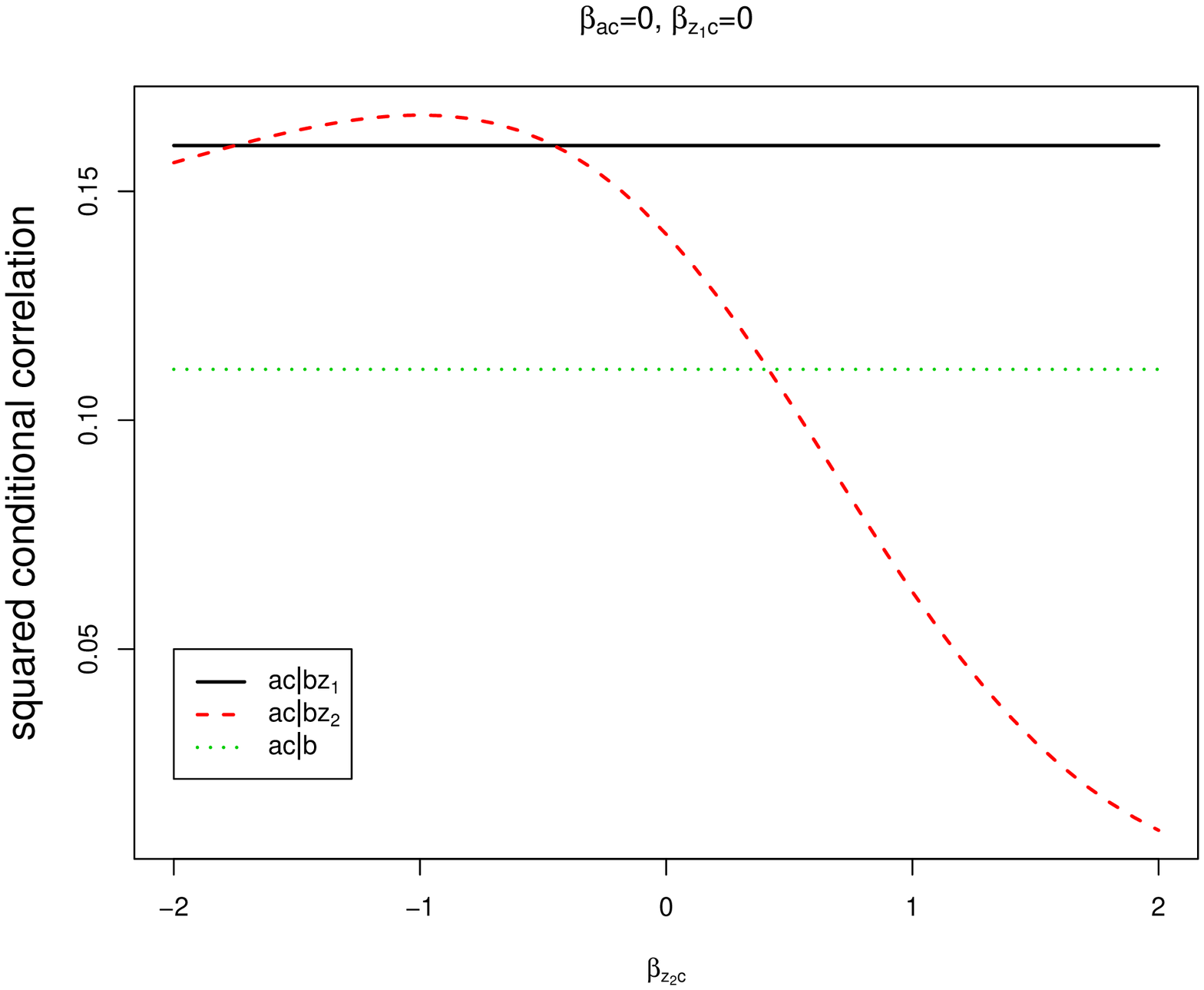}}
}
\caption{The directed acyclic graph described in Example 1. Clearly there is no ordering between $\rho^2_{ac}$, $\rhsq{a}{c}{z_1}$ and $\rhsq{a}{c}{z_2}$.}
\label{fig:2}
\end{figure}

The graphs and the plots used in the counterexamples are described as follows.  In Figures \ref{fig:1} and \ref{fig:2} the graphs with solid edges satisfy the assumptions of Theorems \ref{ap:1b} and \ref{ap:2b} respectively.  We consider the graph with the dashed edges. 
However, excepting one such edge, for all others their corresponding regression coefficients are set to zero.  Each edge implies violation of one conditional independence relationship. 
  
The plots are interpreted as follows. The title of the plots describe which regression coefficients are set to zero. The other regression coefficient is changed and the values of the conditional and unconditional regression coefficients are calculated.  
 
\subsubsection{Figure \ref{fig:1}} The graph with only the solid edges satisfy the conditions of Theorem \ref{ap:1b}.  If an edge between $a$ and $c$ is added, ie. if $\beta_{ac}\ne 0$, but $\beta_{z_1c}=\beta_{z_2a}=0$, $\cind{a}{c}{x}$ no longer holds.  
Figure \ref{fig:c1e1} shows that none of $\rho^2_{ac}$, $\rhsq{a}{c}{z_1}$ and $\rhsq{a}{c}{z_2}$ can be qualitatively compared.  Note that, when $\beta_{ac}=0$ the graph satisfies the condition of Theorem \ref{ap:1b}.  So we get $\rhsq{a}{c}{z_1}\le\rhsq{a}{c}{z_2}\le\rho^2_{ac}$ as predicted. 

If we set $\beta_{ac}=\beta_{z_2a}=0$ and allow $\beta_{z_1c}$ to vary, then for non-zero values of $\beta_{z_1c}$ the condition $\cind{ac}{z_1}{x}$ is violated.  So in figure \ref{fig:c1e2} we see that, the concerned squared partial correlation coefficients are not comparable.

When $\beta_{ac}=\beta_{z_1c}=0$ and $\beta_{z_2a}$ varies, the condition $\cind{z_2}{acx}{z_1}$ is potentially violated. The condition $\cind{z_2}{x}{z_1}$ is not required for Theorem \ref{ap:1b} but for most graphical Markov models $\cind{z_2}{ac}{z_1}$ would imply this condition. Figure \ref{fig:c1e3} shows that the squared correlations cannot be qualitatively compared in this case either.     

The above examples show that none of the conditions of Theorem \ref{ap:1b} can be relaxed further. 

\subsubsection{Figure \ref{fig:2}} In this figure the graph with solid edges satisfy the conditions of Theorem \ref{ap:2b}. If $\beta_{ca}\ne 0$ then the assumption that $a\ind c$ is violated.  As it is evident from the plot in Figure \ref{fig:c2e1} $rhsq{a}{c}{z_1}$, $rhsq{a}{c}{z_2}$ and $\rhsq{a}{c}{x}$ cannot be qualitatively compared.

If $\beta_{z_1c}\ne 0$, $c\not\ind z_1\mid x$ and from Figure \ref{fig:c2e2} it is seen that the squared correlations cannot be qualitatively compared either.

Finally, when $\beta_{z_2c}\ne 0$, $z_2$ becomes conditionally dependent on $c$ given $z_1$.  From Figure \ref{fig:c2e3} we once again conclude that the squared correlations under consideration cannot be qualitatively compared.

The above examples prove that no conditions in Theorem \ref{ap:2b} can be relaxed.

\section{Discussion}
Qualitative comparison may be possible under other sets of conditional independence relations. 
The requirement of a single component $x$ cannot be relaxed. 
The results in Section \ref{sec:canon} are sufficient for postulating path based rules for comparison on polytree models as well. Since the edges on a polytree are directed, these rules are more involved than those for trees \citep{sctsr3}.

Comparison of mutual information with a fixed conditionate holds for any distribution.  In fact, the results with fixed correlates are based on the positive-definiteness of the covariance matrix and extend to non-Gaussian distributions as well. However, inequalities for squared partial correlation would not translate to mutual information for such random variables.  These results may be applicable to causal model selections among non-Gaussian variables (eg. \citet{shimizu:hoyer:etal:2006}). 

It can be shown that, although the comparisons with a fixed conditionate do not hold, but absolute values of partial regression coefficients can be qualitatively compared for fixed correlates under the same conditions \citep{chaudhuri_tan_2010}. 

Rules for signed comparisons of partial correlation and regression coefficients can be developed from these results. Such results might be useful in identifying hidden variables in Factor models \citep{bekker:deleeuw:1987, Drton:strumfels:sullivan:2007, xu:pearl:1989, spirtes:etal:2000} and in recovering population covariance matrix for one-factor models in presence of selection bias \citep{kuroki:cai:2006}.

\appendix

\section{Proofs}\label{sec:proofs}


{\bf Notation:} For two real numbers $a$ and $b$, $a\propto^+ b$ implies that, $\exists$ $M>0$ such that $a=M\cdot b$. 

\begin{prp}\label{prp:covvar}
Suppose $U$, $V$, $W$ are univariate components of a Gaussian random vector with mean $\mu$ and positive definite covariance $\Sigma$. Assume that $\cind{U}{V}{W}$. Then $\sigma_{UV}=\sigma_{UW}\sigma_{WV}/\sigma_{WW}$ and $\sigma_{UU}=\sigma^2_{UW}\sigma_{WW}/\sigma^2_{WW}+E\left[Var\left(U|W\right)\right]$.
\begin{proof} Trivial.\end{proof}
\end{prp}

Suppose $K$ and $K^{\prime}$ are constants and for some $a,c,d \in V$ and $B\subseteq V\setminus\{a,c,d\}$ (where $B$ may be empty) we denote $M_1=\sigma_{cd|B}\left\{\sigma_{ad|B}\sigma_{cc|B}-\sigma_{ac|B}\sigma_{cd|B}\right\}$, \\$M_2=\sigma_{ad|B}\left\{\sigma_{cd|B}\sigma_{aa|B}-\sigma_{ac|B}\sigma_{ad|B}\right\}$, $M_3(\alpha)=[(\alpha-K^{\prime})\sigma_{ac|B}\sigma_{dd|B}-K\cdot \sigma_{ad|B}\sigma_{cd|B}]$ and
\begin{equation}
L(\alpha)=\frac{\{(\alpha-K^{\prime})\rho_{ac|B}-K\rho_{ad|B}\rho_{cd|B}\}^2}{[\{(\alpha-K^{\prime})-K\rho^2_{ad|B}\}\{(\alpha-K^{\prime})-K\rho^2_{cd|B}\}]}\label{eq:l}.
\end{equation}

\begin{lem}\label{lem:initthree}
Suppose $K>0$ and for some $K^{\prime}$ and $\alpha$, $(\alpha-K^{\prime})-K\rho^2_{ad|B}> 0$ and $(\alpha-K^{\prime})-K\rho^2_{cd|B}> 0$. 

Then if $M_1\cdot M_2\ge 0$ :
\begin{enumerate}
\item $\frac{\partial L(\alpha)}{\partial\alpha}=0$ if both $M_1\cdot M_3(\alpha)$ and $M_2\cdot M_3(\alpha)$ are $0$.
\item $\frac{\partial L(\alpha)}{\partial\alpha}$ has the same sign as either $M_1\cdot M_3(\alpha)$ or $M_2\cdot M_3(\alpha)$, whichever is non-zero. 
\end{enumerate}
\begin{proof}
Since the denominator of \eqref{eq:l} is positive then the sign of $\partial L(\alpha)/\partial\alpha$ is the sign of the numerator of
$\partial L(\alpha)/\partial\alpha$. From quotient rule of differentiation and some algebraic manipulation we get :

\begin{align}\label{eq:dl}
\frac{\partial L(\alpha)}{\partial\alpha}\propto^+&K[(\alpha-K^{\prime})\rho_{ac|B}-K\rho_{ad|B}\rho_{cd|B}]\times\nonumber\\
&\left\{[(\alpha-K^{\prime})-K\rho^2_{ad|B}]\rho_{cd|B}[\rho_{ad|B}-\rho_{ac|B}\rho_{cd|B}]\right.\nonumber\\
&+\left.[(\alpha-K^{\prime})-K\rho^2_{cd|B}]\rho_{ad|B}[\rho_{cd|B}-\rho_{ac|B}\rho_{ad|B}]\right\}.
\end{align}
Note that $\rho_{cd|B}[\rho_{ad|B}-\rho_{ac|B}\rho_{cd|B}]\propto^+M_1$, $\rho_{ad|B}[\rho_{cd|B}-\rho_{ac|B}\rho_{ad|B}]\propto^+M_2$ and $[(\alpha-K^{\prime})\rho_{ac|B}-K\rho_{ad|B}\rho_{cd|B}]\propto^+M_3(\alpha)$. 
By substituting these expressions in \eqref{eq:dl} and the positivity $K$, $[(\alpha-K^{\prime})-K\rho^2_{ad|B}]$ and $[(\alpha-K^{\prime})-K\rho^2_{cd|B}]$ the result follows.\hfill$\qed$
\end{proof}
\end{lem}

\noindent\emph{\bf Proof of Theorem \ref{lem:maincomp1}}.
From the assumption $Inf\left(\cind{a}{c^{\prime}}{cZ}\right)=0$. The rest follows from the identity $Inf\left(\cind{a}{c^{\prime}}{cZ}\right)+Inf\left(\cind{a}{c}{Z}\right)=Inf\left(\cind{a}{c}{c^{\prime}Z}\right)+Inf\left(\cind{a}{c^{\prime}}{Z}\right)$.\footnote{The author would like to thank the referee for drawing his attention to this equality which improved the original proof immensely.}\hfill$\qed$

Note that, from \citet{lnenicka_matus_2007}, assumptions on conditional independence and the conditional correlations do not change if we replace $\Sigma$ by $J\Sigma J$, where $J$ is the diagonal matrix with $1/\sqrt{\sigma_{vv}}$, $v\in V$.  Thus, unless otherwise stated, w.l.g we can assume that the diagonal elements of $\Sigma$ are all equal to $1$ and all the off diagonals are in $(-1,1)$.  
That is $\Sigma$ is the correlation matrix of $V$, but with an abuse of notation in what follows below, we still denote the correlation of $a$ and $c$ by $\sigma_{ac}$.

\noindent\emph{\bf Proof of Theorem \ref{ap:1b}.}
Note that by assumption $\sigma_{ac}=\sigma_{ax}\sigma_{cx}$, $\sigma_{az}=\sigma_{ax}\sigma_{xz}$, $\sigma_{cz}=\sigma_{cx}\sigma_{xz}$, $\sigma_{az^{\prime}}=\sigma_{ax}\sigma_{xz}\sigma_{zz^{\prime}}$ and $\sigma_{cz^{\prime}}=\sigma_{cx}\sigma_{xz}\sigma_{zz^{\prime}}$. 

\noindent Part $1$. $\rhsq{a}{c}{z}=\sigma^2_{ac}(1-\sigma^2_{xz})^2/[(1-\sigma^2_{ax}\sigma^2_{xz})(1-\sigma^2_{cx}\sigma^2_{xz})]$. Now since $\sigma^2_{ax}\sigma^2_{xz}\le \sigma^2_{xz}$ and $\sigma^2_{cx}\sigma^2_{xz}\le \sigma^2_{xz}$, $\rhsq{a}{c}{z}\le\rho^2_{ac}$.  

\noindent Part $2$. Assume that $x\ne z^{\prime}$ and consider three non trivial cases as $x=a$, $x=c$ and $x\not\in\{a,c\}$.  Initially assume that $\sigma_{zz^{\prime}}\ne 0$. Since $\cind{ac}{z^{\prime}}{z}$, using Proposition \ref{prp:covvar} and the positive definiteness of the covariance matrix together with $\tau^2_{z}=(1-\sigma^2_{zz^{\prime}})>0$ 
and by denoting $\alpha=1+(\tau_z^2/\sigma^2_{zz^{\prime}})>1$, with $B=\emptyset$, $K^{\prime}=0$, $K=1$ it follows that $\rhsq{a}{c}{z^{\prime}}=L(\alpha)$ for $\alpha\ge 1$ and $\rhsq{a}{c}{z}=L(1)$. Thus in Lemma \ref{lem:initthree} using Cauchy Schwartz inequality and $\alpha\ge 1$ it follows that for $x=a$, $M_1\propto^+\sigma_{cx}$, $M_2=0$ and 
$M_3(\alpha)\propto^+\sigma_{cx}$, for $x=c$, $M_1=0$, $M_2\propto^+\sigma_{ax}$ and $M_3(\alpha)\propto^+\sigma_{ax}$ and for $x\not\in\{a,c\}$, $M_1\propto^+\sigma_{ax}\sigma_{cx}$, $M_2\propto^+\sigma_{ax}\sigma_{cx}$ and $M_3(\alpha)\propto^+\sigma_{cx}\sigma_{ax}$.  
Thus for all cases $\partial L/\partial\alpha\ge 0$ and the result follows.  If $\sigma_{zz^{\prime}}=0$, $z\ind z^{\prime}$ and $z^{\prime}\ind acz$.  Thus $\rhsq{a}{c}{z^{\prime}}=\rho^2_{ac}$.  The rest follows from part $1$.   

For the second inequality notice that, by our assumption $\sigma_{az^{\prime}}=\sigma_{az}\sigma_{zz^{\prime}}=\sigma_{ax}\sigma_{xz}\sigma_{zz^{\prime}}$.  Since we don't assume $\cind{x}{z^{\prime}}{z}$, $\sigma_{xz}\sigma_{zz^{\prime}}$ is not necessarily equal to $\sigma_{xz^{\prime}}$. However, $\rhsq{a}{c}{z^{\prime}}=\sigma^2_{ac}(1-\sigma^2_{xz}\sigma^2_{zz^{\prime}})^2/[(1-\sigma^2_{ax}\sigma^2_{xz}\sigma^2_{zz^{\prime}})(1-\sigma^2_{cx}\sigma^2_{xz}\sigma^2_{zz^{\prime}})]\le\rho^2_{ac}$ in the same way as in part $1$.\hfill$\qed$ 



\noindent\emph{\bf Proof of Theorem \ref{ap:2b}.}
By assumption $\cind{zB}{ac}{x}$ and $a\ind c$.  

Part $1$. It is enough to show that $\sigma^2_{ac|Bz}\ge\sigma^2_{ac|B}$.  Using the above relations in Proposition \ref{prp:covvar} and by denoting $Q_1=\Sigma_{xB}\Sigma^{-1}_{BB}\Sigma_{Bx}$ and $Q_2=\left(\Sigma_{xB},\sigma_{xz}\right)\Sigma^{-1}_{(Bz)(Bz)}\left(\Sigma_{xB},\sigma_{xz}\right)^T$ 
one gets $\sigma_{ac|Bz}=-\sigma_{ax}\sigma_{cx}Q_2$ and $\sigma_{ac|B}=-\sigma_{ax}\sigma_{cx}Q_1$. Now the proof follows by noting that, $\sigma_{aa}-\sigma^2_{ax}Q_1=\sigma_{aa|B}\ge\sigma_{aa|Bz}=\sigma_{aa}-\sigma^2_{ax}Q_2$ implies $Q_2\ge Q_1$.

Part $2$. We initially assume that $\sigma_{zz^{\prime}}\ne 0$.  By defining $\tau^2_{z^{\prime}}=\left(1-\sigma^2_{zz^{\prime}}\right)>0$, $\alpha=\left(1+(\tau^2_{z^{\prime}}/\sigma^2_{zz^{\prime}})\right)$, $K^{\prime}=\Sigma_{zB}\Sigma^{-1}_{BB}\Sigma_{Bz}>0$, $K=(1-K^{\prime})>0$ 
and from the assumption that $\cind{z^{\prime}}{acB}{z}$ it follows that $\rho^2_{ac|Bz^{\prime}}=L(\alpha)$ with $\alpha\ge 1$ and $\rhsq{a}{c}{Bz}=L(1)$. Further using $\cind{ac}{zB}{x}$ one 
can show that $M_1\propto^+\sigma_{cx}\sigma_{ax}\sigma^2_{xz|B}$, $M_2\propto^+\sigma_{cx}\sigma_{ax}\sigma^2_{xz|B}$ and $M_3(\alpha)\propto^+-\sigma_{cx}\sigma_{ax}$.  
Thus from Lemma \ref{lem:initthree} it follows that $\partial L/\partial\alpha\le 0$.  If $\sigma_{zz^{\prime}}=0$, as before $z\ind z^{\prime}$ and $z^{\prime}\ind acB$. Thus $\rhsq{a}{c}{Bz}=\rhsq{a}{c}{B}$. The result follows from part $1$.  

For the first inequality, notice that $\sigma_{ac\mid Bz^{\prime}}=-\sigma_{ax}\sigma_{cx}Q^{\star}_2$ and $\sigma_{aa\mid Bz^{\prime}}=1-\sigma^2_{ax}Q^{\star}_2$, where $Q^{\star}_2=\left(\Sigma_{xB},\sigma_{xz}\sigma_{zz^{\prime}}\right)\Sigma^{-1}_{(Bz^{\prime})(Bz^{\prime})}\left(\Sigma_{xB},\sigma_{xz}\sigma_{zz^{\prime}}\right)^T$.  This implies $\sigma^2_{ac\mid Bz^{\prime}}\ge\sigma^2_{ac\mid B}$ just like part $1$ above.\hfill$\qed$


\noindent\emph{\bf Proof of Theorem \ref{ap:3b}.}
W.l.g. it is enough assume that $x\not\in B$.  Furthermore, note that $\sigma_{aa|B}\ge\sigma_{aa|Bz}$ and $\sigma_{cc|B}\ge\sigma_{cc|Bz}$, thus for part $1$ it is enough to show that under the assumptions $\sigma_{ac|Bz}=m\cdot\sigma_{ac|B}$ for some $m>1$.\\

Part $1$.  Assume that, $a\ind z$ and let (ii) hold, ie. $\cind{cB}{az}{x}$.  Using Proposition \ref{prp:covvar} it follows that 
\begin{align}
\sigma_{ac|Bz}&=\sigma_{ac|B}+\frac{(\Sigma_{aB}\Sigma^{-1}_{BB}\Sigma_{Bz})(\sigma_{cz}-\Sigma_{cB}\Sigma^{-1}_{BB}\Sigma_{Bz})}{\sigma_{zz|B}}=\sigma_{ac|B}+\frac{\sigma_{ax}\sigma^2_{xz}Q_1(\sigma_{cx}-\sigma_{cx}Q_1)}{\sigma_{zz|B}}\nonumber\\
&=\sigma_{ac|B}+\frac{\sigma^2_{xz}Q_1(\sigma_{cx}\sigma_{ax}-\sigma_{cx}\sigma_{ax}Q_1)}{\sigma_{zz|B}}=\sigma_{ac|B}\left(1+\sigma^2_{zx}Q_1\sigma^{-1}_{zz|B}\right).\nonumber 
\end{align}
Thus $\rhsq{a}{c}{B}\le\rhsq{a}{c}{Bz}$.
Under (i) if $c\ind{az}$, $\sigma_{ac}=\sigma_{zc}=0$, $\sigma_{ac|B}=-\Sigma_{aB}\Sigma^{-1}_{BB}\Sigma_{Bc}$ and $\sigma_{cz|B}=-\Sigma_{cB}\Sigma^{-1}_{BB}\Sigma_{Bz}$.  Now if $(i)(a)$ ie. $\cind{az}{B}{x}$ holds:
\begin{align}
\sigma_{ac|Bz}&=\sigma_{ac|B}-\frac{(\Sigma_{aB}\Sigma^{-1}_{BB}\Sigma_{Bz})(\Sigma_{cB}\Sigma^{-1}_{BB}\Sigma_{Bz})}{\sigma_{zz|B}}\label{eq:cov}\\
&=\sigma_{ac|B}-\frac{(\sigma_{ax}\sigma_{xz}Q_1)(\Sigma_{cB}\Sigma^{-1}_{BB}\Sigma_{Bx}\sigma_{xz})}{\sigma_{zz|B}}=\sigma_{ac|B}\left(1+\sigma^2_{zx}Q_1\sigma^{-1}_{zz|B}\right).\nonumber 
\end{align}
Under $(i)(b)$ ie. $\cind{az}{B}{cx}$ notice that from Proposition \ref{prp:covvar}:
\begin{align}
\Sigma_{aB}&=\Sigma_{a(xc)}\Sigma^{-1}_{(xc)(xc)}\Sigma_{(xc)B}=[\sigma_{ax},0]\Sigma^{-1}_{(xc)(xc)}\Sigma_{(xc)B}=\sigma_{ax}[1,0]\Sigma^{-1}_{(xc)(xc)}\Sigma_{(xc)B}=\sigma_{ax}\mathcal{Q}_{cxB}.\nonumber
\end{align}
Here $\mathcal{Q}_{cxB}=[1,0]\Sigma^{-1}_{(xc)(xc)}\Sigma_{(xc)B}$.  Similarly it can be shown that, $\Sigma_{zB}=\sigma_{zx}\mathcal{Q}_{cxB}$ and $\sigma_{ac|B}=-\sigma_{ax}\mathcal{Q}_{cxB}\Sigma^{-1}_{BB}\Sigma_{Bc}$.
Now by substitution in \eqref{eq:cov} above we get:
\begin{align}
\sigma_{ac|Bz}&=\sigma_{ac|B}-\frac{\sigma_{ax}(\mathcal{Q}_{cxB}\Sigma^{-1}_{BB}\mathcal{Q}^T_{cxB})(\Sigma_{cB}\Sigma^{-1}_{BB}\mathcal{Q}^T_{cxB})\sigma^2_{zx}}{\sigma_{zz|B}}\nonumber\\
&=\sigma_{ac|B}-\frac{(\sigma_{zx}\mathcal{Q}_{cxB}\Sigma^{-1}_{BB}\mathcal{Q}^T_{cxB}\sigma_{zx})(\Sigma_{cB}\Sigma^{-1}_{BB}\mathcal{Q}^T_{cxB}\sigma_{ax})}{\sigma_{zz|B}}\nonumber\\
&=\sigma_{ac|B}+\frac{(\Sigma_{zB}\Sigma^{-1}_{BB}\Sigma_{Bz})\sigma_{ac|B}}{\sigma_{zz|B}}=\sigma_{ac|B}\left\{1+(\Sigma_{zB}\Sigma^{-1}_{BB}\Sigma_{Bz})\sigma^{-1}_{zz|B}\right\}.\nonumber
\end{align}
The proofs for $(i)(c)$ and $(i)(d)$ are similar.
 
If $(i)(e)$ ie. $\cind{ac}{B}{x}$ holds, $\sigma_{ac|B}=-\sigma_{ax}\sigma_{cx}Q_1$ and using Proposition \ref{prp:covvar} we get, 
\begin{align}
\sigma_{ac|Bz}=&\sigma_{ac|B}-\frac{(-\Sigma_{aB}\Sigma^{-1}_{BB}\Sigma_{Bz})(-\Sigma_{cB}\Sigma^{-1}_{BB}\Sigma_{Bz})}{\sigma_{zz|B}}=\sigma_{ac|B}-\sigma_{ax}\sigma_{cx}(\Sigma_{xB}\Sigma^{-1}_{BB}\Sigma_{Bz})^2\sigma^{-1}_{zz|B}\nonumber\\
&=\sigma_{ac|B}\left\{1+(\Sigma_{xB}\Sigma^{-1}_{BB}\Sigma_{Bz})^2/(Q_1\sigma^{-1}_{zz|B})\right\}.\nonumber
\end{align}
Under condition $(i)(f)$ notice that, $\Sigma_{aB}=\Sigma_{a(xz)}\Sigma^{-1}_{(xz)(xz)}\Sigma_{(xz)B}=\sigma_{ax}[1,0]\Sigma^{-1}_{(xz)(xz)}\Sigma_{(xz)B}=\sigma_{ax}\mathcal{Q}_{xzB}$.  Similarly, $\Sigma_{cB}=\sigma_{cx}\mathcal{Q}_{xzB}$. Now from \eqref{eq:cov} it follows that:
\[\sigma_{ac|Bz}=\sigma_{ac|B}-\frac{\sigma_{ax}\sigma_{cx}(\mathcal{Q}_{xzB}\Sigma^{-1}_{BB}\Sigma_{Bz})^2}{\sigma_{zz|B}}.
\]
Clearly if at least one of $\sigma_{ax}$,$\sigma_{cx}$, $\mathcal{Q}_{xzB}$ is zero, the results is trivial. Now suppose none of them equal zero.  Then $\mathcal{Q}_{xzB}\Sigma^{-1}_{BB}\mathcal{Q}^T_{xzB}>0$.  Further $\sigma_{ac|B}=-\sigma_{ax}\sigma_{cx}(\mathcal{Q}_{xzB}\Sigma^{-1}_{BB}\mathcal{Q}^T_{xzB})$, which yields
\[\sigma_{ac|Bz}=\sigma_{ac|B}\left\{1+\frac{(\mathcal{Q}^T_{xzB}\Sigma^{-1}_{BB}\Sigma_{Bz})^2}{(\mathcal{Q}_{xzB}\Sigma^{-1}_{BB}\mathcal{Q}^T_{xzB})\sigma_{zz|B}}\right\}.
\] 

Part $2$.  
Suppose $\sigma^2_{z^{\prime}z}>0$. Let $\tau^2_{z^{\prime}}=\left(1- \sigma^2_{z^{\prime}z}\right)>0$, $K^{\prime}=\Sigma_{zB}\Sigma_{BB}^{-1}\Sigma_{Bz}$, $K=(1-K^{\prime})>0$
 and $\alpha=1/\sigma^2_{z^{\prime}z}=\left(1+\tau^2_{z^{\prime}}/\sigma^2_{z^{\prime}z}\right)\ge 1$. Then from $\cind{acB}{z^{\prime}}{z}$, $a\ind zz^{\prime}$ it follows that for both cases $\rho^2_{ac|Bz^{\prime}}=L\left(\alpha\right)$ with $\alpha\ge 1$ and $\rho^2_{ac|Bz}=L(1)$.  Now we consider the four cases in the statement. By denoting $Q_{cx}=\Sigma_{cB}\Sigma^{-1}_{BB}\Sigma_{Bx}$, $Q_{ax}=\Sigma_{aB}\Sigma^{-1}_{BB}\Sigma_{Bx}$ and $\mathcal{Q}_{axB}=[1,0]\Sigma^{-1}_{(xa)(xa)}\Sigma_{(xa)B}$ it follows that:
\begin{equation}
M_1\propto^+M_2\propto^+\begin{cases}
\sigma_{ax}Q_{cx}&\text{if $(i)$, $(a)$}\\
\sigma_{ax}\mathcal{Q}_{cxB}\Sigma^{-1}_{BB}\Sigma_{Bc}&\text{if $(i)$, $(b)$}\\
\sigma_{cx}Q_{ax}&\text{if $(i)$, $(c)$}\\
\sigma_{cx}\mathcal{Q}_{axB}\Sigma^{-1}_{BB}\Sigma_{Ba}&\text{if $(i)$, $(d)$}\\
\sigma_{ax}\sigma_{cx}&\text{if $(i)$, $(e)$}\\
\sigma_{ax}\sigma_{cx}&\text{if $(i)$, $(f)$}\\
-\sigma_{ax}\sigma_{cx}&\text{if $(ii)$}\\
\end{cases},
M_3(\alpha)\propto^+\begin{cases}
-\sigma_{ax}Q_{cx}&\text{if $(i)$, $(a)$}\\
-\sigma_{ax}\mathcal{Q}_{cxB}\Sigma^{-1}_{BB}\Sigma_{Bc}&\text{if $(i)$, $(b)$}\\
-\sigma_{cx}Q_{ax}&\text{if $(i)$, $(c)$}\\
-\sigma_{cx}\mathcal{Q}_{axB}\Sigma^{-1}_{BB}\Sigma_{Ba}&\text{if $(i)$, $(d)$}\\
-\sigma_{ax}\sigma_{cx}&\text{if $(i)$, $(e)$}\\
-\sigma_{ax}\sigma_{cx}&\text{if $(i)$, $(f)$}\\
\sigma_{ax}\sigma_{cx}&\text{if $(ii)$}\\
\end{cases}.\nonumber
\end{equation} 

Thus from Lemma \ref{lem:initthree}, in all cases $\partial L/\partial\alpha\le 0$, which completes the proof. 

If $\sigma_{zz^{\prime}}=0$, then for all cases $\rhsq{a}{c}{Bz^{\prime}}=\rhsq{a}{c}{B}$ and the result follows from Part $1$ as before.

\hfill$\qed$


\noindent\emph{\bf Proof of Corollary \ref{cor:3b}.}
If $B=\emptyset$, under $(i)$ from the assumed independence of $a$, $c$ and $z$, we get $\sigma_{ac}=\sigma_{az}=\sigma_{cz}=0$. The result follows from this.
Under $(ii)$, $\sigma_{cz}\ne 0$ and from Theorem \ref{ap:3b} the result follows.\hfill$\qed$

\noindent\emph{\bf Proof of Theorem \ref{cl:3b2}.}
In this proof we take $\Sigma$ to be the covariance matrix and not the correlation matrix as above. Using condition $\cind{B}{acz}{x}$, denoting $\sigma^2_{xx}Q_4=\Sigma_{xB}\Sigma^{-1}_{BB}\Sigma_{Bx}$, $T=\sigma_{zz}/\left(\sigma_{zz}-\sigma^2_{xz}Q_4\right)$ ($T>0$) and from Proposition \ref{prp:covvar} and some simplification we get 

\begin{equation*}
\frac{\rhsq{a}{c}{Bz}}{\rhsq{a}{c}{x}}=\frac{\left(\sigma_{aa}\sigma_{xx}Q_4T-\sigma^2_{ax}Q_4T\right)\left(\sigma_{cc}\sigma_{xx}Q_4T-\sigma^2_{cx}Q_4T\right)}{\left(\sigma_{aa}-\sigma^2_{ax}Q_4T\right)\left(\sigma_{cc}-\sigma^2_{cx}Q_4T\right)}.
\end{equation*}
Thus $\rhsq{a}{c}{Bz}\ge\rhsq{a}{c}{x}$ iff $\sigma_{xx}Q_4T\ge 1$ iff $\left(\sigma_{xx}+\sigma^2_{xz}/\sigma_{zz}\right)Q_4\ge 1$.  The equivalent expression follows as: 
\begin{equation}
\frac{\sigma^2_{xz}}{\sigma_{zz}\sigma^2_{xx}}\Sigma_{xB}\Sigma^{-1}_{BB}\Sigma_{Bx}\ge\frac{\sigma_{xx|B}}{\sigma_{xx}}\Leftrightarrow\frac{1}{\sigma_{zz}}\Sigma_{zB}\Sigma^{-1}_{BB}\Sigma_{Bz}\ge\frac{\sigma_{xx|B}}{\sigma_{xx}}\Leftrightarrow\frac{\sigma_{xx}-\sigma_{xx|B}}{\sigma_{xx}}\ge\frac{\sigma_{zz|B}}{\sigma_{zz}}.\nonumber
\end{equation}\hfill$\qed$


\noindent\emph{\bf Proof of Lemma \ref{lem:ugcond}.}
We need to show that if $T_1$, $T_2$ and $T_3$, are disjoint subsets of $Z^c$, then $T_1$ is connected to $T_2$ given $T_3$ in $G_{Z^c}$ iff $T_1$ is connected to $T_2$ given $T_3\cup Z$ in $G$.

($\Rightarrow$)  Suppose $T_1$ is connected to $T_2$ given $T_3$ in $G_{Z^c}$.  So there are $t_1\in T_1$ and $t_2\in T_2$ and the path \path{t_1}{t_2} such that $\epath{}{}\cap T_3=\emptyset$.  Clearly \path{t_1}{t_2} is in $G$ and $\epath{t_1}{t_2}\cap Z=\emptyset$. 
So $\epath{t_1}{t_2}\cap \{T_3\cup Z\}=\emptyset$.  This shows $T_1$ is connected to $T_2$ given $T_3\cup Z$ in $G$.

($\Leftarrow$) Suppose $T_1$ is connected to $T_2$ given $T_3\cup Z$ in $G$.  So there is $t_1\in T_1$ and $t_2\in T_2$ and the path \path{t_1}{t_2}, such that $\epath{t_1}{t_2}\cap \{T_3\cup Z\}=\emptyset$. 
So $\epath{t_1}{t_2}\cap Z=\emptyset$ and \path{t_1}{t_2}$\subseteq Z^c$.  Clearly in $G_{Z^c}$, $\epath{t_1}{t_2}\cap T_3=\emptyset$.  This shows $T_1$ is connected to $T_2$ given $T_3$ in $G_{Z^c}$.\hfill$\square$ 

\noindent\emph{\bf Proof of Theorem \ref{thm:cond}.}
From the structure of $G$ and since $c\in\epath{a}{c^{\prime}}$, it easily follows that $c^{\prime}$ is separated from $a$ given $c$ and $Z$. The result follows from Theorem \ref{lem:maincomp1}.\hfill$\square$
 
\noindent\emph{\bf Proof of Theorem \ref{thm:ugcorr}.}
For notational convenience we express the squared partial correlations as functions of the covariance matrix $\Sigma$. We need to show that $\rhsq{a}{c}{Z_1}\left(\Sigma\right)\le \rhsq{a}{c}{Z_2}\left(\Sigma\right)$.  
W.l.g. we assume that for $i=1,2$ there is no $z_i\in Z_i$ such that $\cind{ac}{Z_i\setminus\{z_i\}}{z_i}$.  We consider several cases below:

\noindent\emph{Case $1$.} If $Z_1\cap\epath{a}{c}\ne\emptyset$, then $\cind{a}{c}{Z_1}$, $\rhsq{a}{c}{Z_1}=0$ and the result is trivial.

We initially assume that $Z_1$ separates $Z_2$ from $a$ and $c$. This implies that for each $z_2\in Z_2$ there is a $z^a_1\in Z_1$ and $z^c_1\in Z_1$ such that $z^a_1\in \epath{a}{z_2}$ and $z^c_1\in \epath{c}{z_2}$.

\noindent\emph{Case $2$.} If $Z_2\cap\epath{a}{c}\ne\emptyset$, then $z^a_1\in\epath{a}{z_2}\subseteq\epath{a}{c}$. This implies that $\cind{a}{c}{Z_1}$ and $\rhsq{a}{c}{Z_1}=0$.\\
\noindent\emph{Case $3$.} Now let $\left(Z_1\cup Z_2\right)\cap\epath{a}{c}=\emptyset$. Suppose $Z_1=\left\{z_{11}\right.$, $z_{12}$, $\ldots$, $\left.z_{1n_1}\right\}$ and $Z_2=\left\{z_{21}\right.$, $z_{22}$, $\ldots$, $\left.z_{2n_2}\right\}$. 
Suppose $x_i=\epath{a}{z_{1i}}\cap\epath{c}{z_{1i}}\cap\epath{a}{c}$. Since $G$ is a tree $x_i$ is unique for $z_i$.  Also suppose that $N^{}_i=\left\{z_{2i}\in Z_2 ~:~z_{1i}\in\epath{a}{z_{2i}}\cap\epath{c}{z_{2i}}\right\}$. 
Again from the structure of $G$ it is clear that $N^{}_i$ are disjoint and $Z_2=\cup^{n_1}_{i=1} N^{}_i$. We don't exclude the possibility that $N^{}_i$ may be $\emptyset$ for some $i$.  Using \eqref{eq:factmain} we can write:
\begin{equation}\label{eq:fact}
\frac{\rhsq{a}{c}{Z_1}}{\rhsq{a}{c}{Z_2}}=\prod^{n_1}_{i=1}\frac{\rhsq{a}{c}{z_{11}\ldots z_{1(i-1)}z_{1i}N^{}_{i+1}\ldots N^{}_{n_1}}}{\rhsq{a}{c}{z_{11}\ldots z_{1(i-1)}N^{}_{i}N^{}_{i+1}\ldots N^{}_{n_1}}}.
\end{equation}

It is sufficient to show that each factor in the product \eqref{eq:fact} is bounded by $1$. Consider the $ith$ factor,
\begin{equation*}
f_i=\frac{\rhsq{a}{c}{z_{11}\ldots z_{1(i-1)}z_{1i}N^{}_{i+1}\ldots N^{}_{n_1}}}{\rhsq{a}{c}{z_{11}\ldots z_{1(i-1)}N^{}_{i}N^{}_{i+1}\ldots N^{}_{n_1}}}.
\end{equation*}
Notice that the factor $f_i$ depends only on the subgraph $G_{V_i}$ of $G$ defined by the vertex set:
\begin{equation*}
V_i=\left\{\bigcup^{i-1}_{j=1}\left(\epath{a}{z_{1j}}\cup\epath{c}{z_{1j}}\right)\right\}\bigcup\left\{\bigcup^{n_1}_{j=i}\bigcup_{z^{(j)}_{2k}\in N^{}_j}\left(\epath{a}{z^{(j)}_{2k}}\cup\epath{c}{z^{(j)}_{2k}}\right)\right\}.
\end{equation*}

It is clear that, $G_{V_i}$ is a tree. Let us denote $B_i=\{z_{11},\ldots,z_{1(i-1)}\}\cup\left(\cup^{n_1}_{j=i+1}N^{}_j\right)$ and $B^c_i=V_i\setminus B_i$.   

Now from the structure of $G_{V_i}$ we note that (i) $x_i\in\epath{a}{c}$ so $\cind{a}{c}{x_iB_i}$, (ii) $x_i\in\epath{a}{z_{1i}}$ and $x_i\in\epath{c}{z_{1i}}$ implying $\cind{a}{z_{1i}}{x_i,B_i}$ and 
(iii) $z_{1i}\in\bigcup_{z^{(i)}_{2k}\in N^{}_i}\left(\epath{a}{z^{(i)}_{2k}}\cap\epath{c}{z^{(i)}_{2k}}\right)$ it follows that $\cind{ac}{N^{}_i}{z_{1i}B_i}$.

From Lemma \ref{lem:ugcond} it follows that the triples $\trip{a}{c}{x_i}$, $\trip{ac}{z_{1i}}{x_i}$ and $\trip{ac}{N_i}{z_{1i}}$ are in $\cm{\indrels{G}}{B_i}{\empty}=\mathfrak{I}\left(G_{B_i^c}\right)$.  
It is obvious that, 
\begin{equation*}
\frac{\rhsq{a}{c}{B_iz_{1i}}\left(\Sigma\right)}{\rhsq{a}{c}{B_iN_i}\left(\Sigma\right)}=\frac{\rhsq{a}{c}{z_{1i}}\left(\Sigma_{B^c_iB^c_i|B_i}\right)}{\rhsq{a}{c}{N_i}\left(\Sigma_{B^c_iB^c_i|B_i}\right)}.
\end{equation*}


Now consider the following sub-cases:
\begin{list}{}{}
\item[a.] If $N^{}_i=\emptyset$ or $N^{(1)}_i=z_{2i}$, from the Theorem \ref{ap:1b} it follows that $\rhsq{a}{c}{z_{1i}}\left(\Sigma_{B^c_iB^c_i|B_i}\right)\le \rhsq{a}{c}{N_i}\left(\Sigma_{B^c_iB^c_i|B_i}\right)$.
\item[b.] If $N^{}_i=\{z_{21},\ldots,z_{2m_i}\}$, then using $\cind{ac}{N^{}_i}{z_{1i}}$, we can write:
\begin{equation*}
f_i=\frac{\rhsq{a}{c}{z_{1i}z_{22}\ldots z_{2m_i}}\left(\Sigma_{B^c_iB^c_i|B_i}\right)}{\rhsq{a}{c}{z_{21}z_{22}\ldots z_{2m_i}}\left(\Sigma_{B^c_iB^c_i|B_i}\right)}.
\end{equation*}
By following the same argument as above and conditioning on $\{z_{22},\ldots,z_{2m_i}\}$ it follows that $f_i\le 1$. 
\end{list}

Now suppose that there is a $Z^{\prime}_2\subseteq Z_2$ s.t. $Z^{\prime}_2$ is not separated from $a$ and $c$ by $Z_1$, but because of the choice of parameters both $\rhsq{a}{Z^{\prime}_2}{Z_1}=\rhsq{c}{Z^{\prime}_2}{Z_1}=0$. 

It can be shown that $\rhsq{a}{c}{(Z^{\prime}_2\cup Z_1)}=\rhsq{a}{c}{Z_1}$.  So if $Z^{\prime}_2\cap\epath{a}{c}\ne\emptyset$ then $\rhsq{a}{c}{Z_2}=\rhsq{a}{c}{(Z^{\prime}_2\cup Z_1)}=\rhsq{a}{c}{Z_1}=0$.  On the other hand if $Z^{\prime}_2\cap\epath{a}{c}=\emptyset$ we can write:
\begin{equation}\label{eq:nosep}
\frac{\rhsq{a}{c}{Z_1}}{\rhsq{a}{c}{Z_2}}=\frac{\rhsq{a}{c}{(Z^{\prime}_2\cup Z_1)}}{\rhsq{a}{c}{(Z^{\prime}_2\cup \{Z_2\setminus Z^{\prime}_2\})}}.
\end{equation} 
The fact that the ratio in \eqref{eq:nosep} is less than $1$ follows from the first part mutatis mutandis. \hfill$\square$\\

\noindent\emph{\bf Proof of Corollary \ref{cor:ugcorr}.} The assumptions imply that $Z_1$ separates $Z_2$ from $a$ and $c$.  This is exactly Case $3$. in the previous proof.
\hfill$\square$

\noindent\emph{\bf Proof of Theorem \ref{thm:cmplt}.} We parametrise the Choleski decomposition $\Lambda=BB^T$.  

Suppose $z_1\in Z_1$ and $z_2\in Z_2$ such that $ac\not\ind z_1|Z_2$ and $ac\not\ind z_2|Z_1$.  Let $\epath{a}{c}=\left\{a=v_1,v_2,\ldots,v_d=c\right\}$, $\epath{a}{c}\cap\epath{a}{z_1}\cap\epath{c}{z_1}=v_i$, $\epath{a}{c}\cap\epath{a}{z_2}\cap\epath{c}{z_2}=v_j$, $i,j\in\{1,2,\ldots,d\}$.  
Further let $\epath{v_i}{z_1}=\{v_i,x_1,\ldots,x_{d_1}=z_1\}$ and $\epath{v_i}{z_2}=\{v_j,y_1,\ldots,y_{d_2}=z_2\}$.  If $i=j$ it is possible that \path{v_i}{z_1} and \path{v_i}{z_2} intersect at more than one vertex.  However, it does not change the proof, so w.l.g. we assume that $i\ne j$.  Suppose 
\begin{align}
V_I&=\epath{a}{c}\cup\epath{v_i}{z_1}\cup\epath{v_j}{z_2}\nonumber\\
E_I&=\{(v_2,v_1),\ldots,(v_d,v_{d-1}),(x_1,v_i),\ldots,(z_1,x_{d_1-1}),(y_1,v_j),\ldots,(z_2,y_{d_2-1})\}.\nonumber 
\end{align}
  
We list the variables in $\Sigma$ as $\epath{a}{c},\epath{v_1}{z_1},\epath{v_2}{z_2},V\setminus V_I$, where the vertices in $V\setminus V_I$ can be arranged in an arbitrary fashion. The matrix $B$ inherits the same arrangement. 

The matrix $B$ is given by, $B_{kl}=1$, $\text{\{if $k=l$\}}$, $B_{kl}=-1$, $\text{\{if $(k,l)\in E_I$\}}$, $B_{kl}=-b_1$, $\text{\{if $(k,l)=(z_1,x_{d_1-1})$\}}$, $B_{kl}=-b_2$, $\text{\{if $(k,l)=(z_2,x_{d_2-1})$\}}$, $B_{kl}=0$, $\text{\{otherwise\}}$.


It can be shown that the resulting $\Lambda$ is a n.n.d. matrix for all values of $b_1$ and $b_2$ and will represent all the conditional independence relations on the tree under consideration.  


Now choose $b_1=0$. This implies $\rhsq{a}{c}{Z_1}=\rho^2_{ac}\ge\rhsq{a}{c}{z_2}=\rhsq{a}{c}{Z_2}$.  The opposite happens if $b_2=0$.  
This completes the proof.\hfill$\square$

\noindent\emph{\bf Proof of Corollary \ref{cor:final}.} The result is trivial if $\epath{a}{c^{\prime}}\cap Z^{\prime}\ne\emptyset$. Furthermore, by assumption if $Z$ intersects \path{a}{c}, so does $Z^{\prime}$. The non-trivial case can be shown by applying Theorem \ref{thm:cond} and Corollary \ref{cor:ugcorr} respectively on the factors below:
\begin{equation*}
\frac{\rhsq{a}{c^{\prime}}{Z^{\prime}}}{\rhsq{a}{c}{Z}}=\frac{\rhsq{a}{c^{\prime}}{Z^{\prime}}}{\rhsq{a}{c}{Z^{\prime}}}\frac{\rhsq{a}{c}{Z^{\prime}}}{\rhsq{a}{c}{Z}}.
\end{equation*}
\hfill$\square$
~\hfill\\
To prove Theorem \ref{thm:modsel} we need the following definitions from the literature of directed acyclic graphs.

\begin{defn}\label{defn:coll} A vertex $v$ on a path \path{x}{y} in a polytree is a collider on the path if there are vertices $v_1$ and $v_2$ on \path{x}{y} such that the edges $v_1\rightarrow v$ and $v_2\rightarrow v$ exist.  A vertex on a path \path{x}{y} in a polytree is a non-collider on the path if it is not a collider on \path{x}{y}.
\end{defn}

\begin{defn}[d-connection]\label{defn:dconn} A path \path{x}{y} between $x$ and $y$ in a DAG is said to be d-connecting given a set $Z$ (possibly empty) if $1.$ every non-collider on \path{x}{y} is not in $Z$ and $2.$ every collider on \path{x}{y} is in $an(Z)$. Here $an(Z)=\cup_{z\in Z}an(z)$.

If there is no path d-connecting $x$ and $y$ given $Z$, then $x$ and $y$ are said to be d-separated given $Z$. 
\end{defn}
\begin{defn}For disjoint sets $X$, $Y$, $Z$, where $Z$ may be empty, $X$ and $Y$ are d-separated given $Z$, if for every pair $x$, $y$, with $x\in X$ and $y\in Y$, $x$ and $y$ are d-separated given $Z$.
\end{defn}
\begin{defn} We say a density $f$ factors according to a DAG, if for three disjoint sets $X$, $Y$ and $Z$, $\cind{X}{Y}{Z}$ according to $f$ whenever $X$ is d-separated from $Y$ given $Z$. 
\end{defn}
\begin{figure}[ht]
\begin{center}
\input{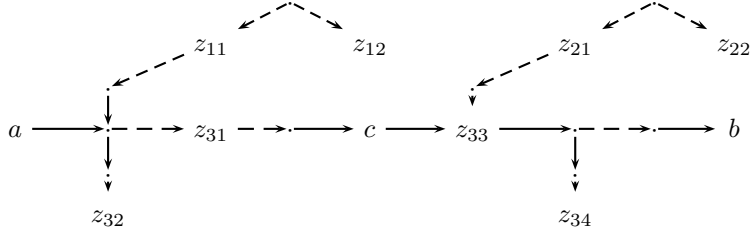}
\caption{Example of a polytree discussed in Theorem \ref{thm:modsel}. Vertices $z_{11}$ and $z_{12}$ are relevant to Case $(i)$, $z_{21}$ and $z_{22}$ are relevant to Case $(ii)$ below. The vertices $z_{3k}$, for $k=1,\ldots,4$ corresponds to Case $(iii)$ in the proof below.}
\end{center}
\end{figure}

\noindent \emph{\bf Proof of Theorem \ref{thm:modsel}.} First of all note that, since $\rhsq{a}{c}{bz}\ne \rhsq{a}{c}{b}$, $ac\not\ind z\mid b$. Further, since $a\ne c\ne b$, $a\in an(c)$ and $c\in an(b)$, there are no colliders on \path{a}{b}.  We first show that $z\ind ac$ iff $z\ind c$. Clearly, $z\ind ac$ implies $z\ind c$. 
To show the converse first note that, since the graph is a polytree, if $z\ind c$ there is at least one collider $v$ on the unique path \path{c}{z} between $c$ and $z$.  Clearly, $v$ cannot be on \path{a}{c}, otherwise it will be a collider on \path{a}{c}.  
However, by construction $\epath{c}{z}\setminus\epath{a}{c}=\left(\epath{a}{z}\cap\epath{c}{z}\right)\setminus\epath{a}{c}$. So if $v$ is not on \path{a}{c}, $v$ would be a collider on \path{a}{z} as well. Thus, using the assumption that the graph is a polytree, $a\ind z$ and our claim follows.  

Similar argument shows if $z\ind b$ iff $z\ind acb$. So, $\rhsq{a}{c}{bz}\ne \rhsq{a}{c}{b}$ implies that $z\not\ind b$.  
So only the following three cases, $(i)$ $a\ind z$ and $c\not\ind z$, $(ii)$ $ac\ind z$ (ie. $z\ind c$) and $(iii)$ $a\not\ind z$ and $c\not\ind z$ are possible.  We first consider the if parts:

\noindent Case $(i)$ We show that there is a vertex $v_1$ such that $\cind{az}{cb}{v_1}$.  $a\ind z$ implies there is at least one collider $v_1$ on \path{a}{z}, $a\ne z\ne v_1$.  Again by construction $\epath{c}{v_1}\setminus\epath{a}{c}=\left(\epath{a}{v_1}\cap\epath{c}{v_1}\right)\setminus\epath{a}{c}$. Thus, if $v_1\not\in\epath{a}{c}$, $v_1$ is a collider on \path{c}{z} as well, which would imply $c\ind z$. Thus $v_1\in\epath{a}{c}$.  Clearly, $v_1$ cannot be a collider on \path{a}{c}.
Thus $v_1$ is the only collider on \path{a}{z} and it is not a collider on \path{a}{c} and \path{c}{z}. Thus, from the definition of d-separation  it follows that$\cind{az}{cb}{v_1}$. From Theorem \ref{ap:3b} $(ii)$ it follows that $\rhsq{a}{c}{bz}>\rhsq{a}{c}{b}$.
      
\noindent Case $(ii)$ We show that $\cind{a}{z}{cb}$ and apply Theorem \ref{ap:2b} with $x=c$. Since by assumption $c\ind z$ and $b\not\ind z$, as in Case $(i)$ above there is a vertex $v_2\in \epath{c}{b}$ such that $v_2$ is a collider on \path{c}{z} but not a collider on \path{b}{z}. Note that, $v_2\ne c$ or $v_2\ne z$. 
Thus $c$ is a non-collider on both \path{a}{z} and \path{a}{b} and $c$ d-separates $a$ from $\{b,z\}$. This implies $\cind{a}{bz}{c}$, which in turn gives $\cind{a}{z}{cb}$.  Now from Theorem \ref{ap:2b} we get $\rhsq{a}{c}{bz}<\rhsq{a}{c}{b}$.

\noindent Case $(iii)$ Since $a\not\ind z$, it follows that $c\not\ind z$ and $b\not\ind z$. This implies there is no collider on \path{a}{z}, \path{c}{z} and \path{b}{z}.  
Let $v_3=\epath{a}{z}\cap\epath{c}{z}\cap\epath{b}{z}\cap\epath{a}{b}$.  Clearly, $v_3$ is a non-collider on all these paths. So, it follows that $\cind{acb}{z}{v_3}$ \citep[page 29]{lau}.  This implies $\cind{ac}{z}{bv_3}$.  
Further, if $v_3\in\epath{a}{c}$, $\cind{a}{cb}{v_3}$ and $\cind{a}{c}{bv_3}$. 
It is possible that $z=v_3$.
Now if $v_3\in\epath{a}{c}$, Theorem \ref{ap:2b} with $x=v_3$ imply $\rhsq{a}{c}{bz}<\rhsq{a}{c}{b}$. Note that in this case if $v_3=z$, $\rhsq{a}{c}{bz}=0$.  If $v_3\not\in\epath{a}{c}$, we consider two cases. 
\noindent Case (a) $z=v_3\in\epath{c}{b}$. Clearly $\cind{ac}{b}{z}$. Now using Theorem \ref{ap:2b} we get $\rhsq{a}{c}{bz}=\rhsq{a}{c}{z}<\rhsq{a}{c}{b}$. 
\noindent Case (b) When $v_3\not\in\epath{c}{b}$ use Theorem \ref{ap:2b} on conditional covariance given $b$ with $x=c$ to get $\rhsq{a}{c}{bz}<\rhsq{a}{c}{b}$.

The only if parts follow from the if part and the fact that the above three are only possible cases under our assumptions.
  
\hfill{$\square$}
\section{Mixed ancestral graphs}\label{sec:mags}
In this supplement we briefly discuss mixed ancestral graphs. Our discussion closely follows \citet{thomas1}. We also refer to the same text for a more detailed treatment of the class of these graphs.

A graph $G$ is an ordered pair $(V,E)$ where $V$ is a set of vertices and $E$ is a set of edges.

A mixed graph is a graph containing three types of edges, undirected ($\uned$), directed ($\ded$) and bidirected ($\bded$).  The following terminology is used to describe relations between variables in such a graph:
\begin{enumerate}
\item If $\alpha\uned\beta$ in $G$, then $\alpha$ is a neighbour of $\beta$ and $\alpha\in ne(\beta)$.
\item If $\alpha\ded\beta$ in $G$, then $\alpha$ is a parent of $\beta$ and $\alpha\in pa(\beta)$.
\item If $\beta\ded\alpha$ in $G$, then $\alpha$ is a child of $\beta$ and $\alpha\in ch(\beta)$.
\item If $\alpha\bded\beta$ in $G$, then $\alpha$ is a spouse of $\beta$ and $\alpha\in sp(\beta)$.
\end{enumerate}
\begin{defn} A vertex $\alpha$ is said to be an ancestor of a vertex $\beta$ if either there is a directed path $\alpha\ded\cdots\ded\beta$ from $\alpha$ to $\beta$, or $\alpha=\beta$. Further, for $X\subseteq V$ its ancestor set is defined as:
\[an(X)=\{\alpha~:~\alpha\text{ is an ancestor of $\beta$ for some $\beta\in X$}\}. 
\] 
\end{defn} 
\begin{defn} A vertex $\alpha$ is said to be anterior to a vertex $\beta$ if there is a path \path{\alpha}{\beta} on which every edge is either of the form $\gamma\uned\delta$, or $\gamma\ded\delta$ with $\delta$ between $\gamma$ and $\beta$, or $\alpha=\beta$; that is, 
there are no edges $\gamma\bded\delta$ and there are no edges $\delta\ded\gamma$ pointing toward $\alpha$.  Further, for $X\subseteq V$ its anterior set is defined as:
\[ant(X)=\{\alpha~:~\alpha\text{ is an anterior to $\beta$ for some $\beta\in X$}\}. 
\]  
\end{defn} 
\begin{defn} An ancestral graph $G$ is a mixed graph in which the following conditions hold for all vertices $\alpha$ in $G$:
\begin{enumerate} 
\item $\alpha\not\in ant\left(pa(\alpha)\cup sp(\alpha)\right)$ and
\item if $ne(\alpha)\ne\emptyset$ then $pa(\alpha)\cup sp(\alpha)=\emptyset$.
\end{enumerate}
\end{defn}

The d-separation criterion for DAGs can be extended to m-separation criterion for mixed ancestral graphs.

A non-endpoint vertex $\zeta$ on a path is a collider on the path if the edges preceding and succeeding $\zeta$ on the path have an arrowhead at $\zeta$, ie., $\rightarrow\zeta\leftarrow$, $\bded\zeta\bded$, $\bded\zeta\leftarrow$, $\ded\zeta\bded$.  
A non-endpoint vertex $\zeta$ on a path which is not a collider is a noncollider on the path. 

A path between vertices $\alpha$ and $\beta$ in an ancestral graph $G$ is said to be m-connecting given a set $Z$ (possibly empty), with $\alpha$, $\beta\not\in Z$ if:

\begin{enumerate}
\item every noncollider on the path is not in $Z$, and
\item every collider on the path is in the $ant(Z)$. 
\end{enumerate}

If there is no path m-connecting $\alpha$ and $\beta$ given $Z$, then $\alpha$ and $\beta$ are said to be m-separated given $Z$.  Non empty sets $X$ and $Y$ are m-separated given Z, if for every pair $\alpha$, $\beta$ with $\alpha\in X$ and $\beta\in Y$, $\alpha$ and $\beta$ are m-separated given $Z$ ($X$, $Y$ and $Z$ are disjoint sets).

A distribution $F$ is said to satisfy the conditional independence relations represented by a mixed ancestral graph if for disjoint subsets $X$, $Y$ and $Z$, $\cind{X}{Y}{Z}$ according to $F$ whenever $X$ is m-separated from $Y$ given $Z$.

\begin{figure}[t]
\begin{center}
\subfigure[\label{fig:magb}]{\input{mag.tex}}\hspace{1in}
\subfigure[\label{fig:ap3b1cab}]{\input{figap3b1ca.tex}}
\end{center}
\caption{}
\end{figure} 
\subsection{Examples of mixed ancestral graphs in the main text}
\begin{example} Consider the Mixed ancestral graph in Figure \ref{fig:magb}.  There are more than one paths connecting $a$ and $c$. Each of them has a collider on it. As for example, $y_1$ is a collider on the path $\{a,y_1,c\}$. So $a$ is m-separated from $c$ given $\emptyset$.  Thus $a\ind c$.  
Further note that, $x_2$ is a noncollider on each path connecting $\{a,c\}$ and $z$. Thus, $\cind{ac}{z}{x_2}$. Similarly, $\cind{ac}{z^{\prime}}{z}$.
\end{example}
\begin{example}
Now we consider the graph in Figure \ref{fig:ap3b1cab}.  Clearly $a\ind c$. $x$ is a collider on the paths $\{a,x,z\}$ and $\{c,x,z\}$.  Further, $b$ is a collider on the paths $\{a,x,b,z\}$ and $\{c,x,b,z\}$.  So $b$ and $x$ m-separates $a$ and $c$ from $z$ given $\emptyset$.  So $ac\ind z$.  
Now note that, $x$ is a noncollider on the paths $\{a,x,b\}$ and $\{c,x,b\}$.  Also $z$ is a collider on the paths $\{a,x,z,b\}$ and $\{c,x,z,b\}$.  This implies $\{a,c\}$ is m-separated from $b$ given $x$, but not given $zx$.
\end{example}

\noindent {\large\bf Acknowledgement}
The author would like to thank Michael Perlman, Thomas Richardson, Mathias Drton, Antar Bandyopadhyay, the referees and the associate editor for their useful comments and suggestions during the preparation of this article.

\end{document}